\documentclass[12pt]{amsart}

\usepackage{amsmath,amsthm,amsfonts,xcolor}
\usepackage{tikz}
\usetikzlibrary{shapes}
\usetikzlibrary{arrows}
\usetikzlibrary{decorations.markings}

\tikzstyle heightone=[scale=.7,shift={(0,-.3)}]
\tikzstyle heightones=[scale=.8,xscale=.35,shift={(0,.1)}]
\tikzstyle heightoneonehalf=[scale=.9,shift={(0,-.2)}]
\tikzstyle heighttwo=[scale=.9,shift={(0,-.4)}]
\tikzstyle heighttwos=[scale=.5,xscale=.6,shift={(0,-.1)}]
\tikzstyle heightthree=[scale=.6,shift={(0,-.9)}]
\tikzstyle heightthrees=[scale=.4,xscale=.7,shift={(0,-.2)}]

\tikzstyle arrowstyle=[blue,semitransparent,scale=2]

\tikzstyle basiclabel=[draw=none,fill=none,shape=rectangle,inner sep=2pt,scale=.8]
\tikzstyle leftlabel=[basiclabel,anchor=east]
\tikzstyle rightlabel=[basiclabel,anchor=west]
\tikzstyle bottomlabel=[basiclabel,anchor=north]
\tikzstyle toplabel=[basiclabel,anchor=south]

\tikzstyle vertex=[circle,draw,fill=black,inner sep=1pt]
\tikzstyle ciliation=[circle,draw=none,fill=red,inner sep=1pt,semitransparent]
\tikzstyle ciliatednode=[vertex,pin={[pin distance=1mm,pin edge={semitransparent,red},ciliation]#1:{}}]

\tikzstyle vector=[black,thick,rectangle,draw=gray!50!yellow,top color=yellow!30,bottom color=black!10,scale=.8,inner sep=2pt]
\tikzstyle small vector=[vector,scale=.8]
\tikzstyle plain vector=[rectangle,draw=none,fill=white,scale=.7]

\tikzstyle my signal=[black,thick,signal,signal pointer angle=120,draw=blue!50,top color=blue!20,bottom color=black!10,scale=.8,inner sep=2pt]
\tikzstyle matrix=[my signal,signal from=south,signal to=north]
\tikzstyle reverse matrix=[my signal,signal from=north,signal to=south]

\tikzstyle small matrix=[matrix,scale=.7]
\tikzstyle reverse small matrix=[reverse matrix,scale=.7]
\tikzstyle matrix on edge=[small matrix,sloped,rotate=-90]
\tikzstyle reverse matrix on edge=[small matrix,sloped,rotate=90]

\tikzstyle with matrix=[postaction={decorate,decoration={markings,
    mark=at position .5 with {\node[matrix]{#1};}}}]
\tikzstyle with small matrix=[postaction={decorate,decoration={markings,
    mark=at position .5 with {\node[matrix on edge]{#1};}}}]
\tikzstyle with reverse small matrix=[postaction={decorate,decoration={markings,
    mark=at position .5 with {\node[reverse matrix on edge]{#1};}}}]

\tikzstyle trivalent=[very thick]
\tikzstyle dotdotdot=[decorate,decoration={markings,
    mark=at position .3 with{\node{.};},
    mark=at position .5 with {\node{.};},
    mark=at position .7 with {\node{.};}}]

\tikzstyle wavyup=[out=90,in=-90]
\tikzstyle wavydown=[out=-90,in=90]

\tikzstyle permutation=[rectangle,fill=black,draw=black]
\tikzstyle symmetrizer=[rectangle,fill=black,draw=black]
\tikzstyle antisymmetrizer=[rectangle,fill=gray!10,draw=black]

\tikzstyle{every picture}=[semithick,baseline=0pt,heightone,label distance=-1mm]

\title{Trace Diagrams, Signed Graph Colorings, and Matrix Minors}
\author{Steven Morse and Elisha Peterson}
\address{Department of Mathematical Sciences, United States Military Academy, West Point, NY 10996-1905, \texttt{elisha.peterson@usma.edu}}

\def\C{\mathbb{C}}
\def\F{\mathbb{F}}
\def\bfa{\mathbf{a}}
\def\bfu{\mathbf{u}}
\def\bfv{\mathbf{v}}
\def\bfw{\mathbf{w}}

\def\tensor{\otimes}

\def\bse{\mathbf{\hat e}}
\def\sgn{\mathrm{sgn}}

\def\tr{\mathsf{tr}}
\def\adj{\mathsf{adj}}
\def\det{\mathsf{det}}
\def\int{\mathsf{int}}

\def\tdsign#1{(-1)^{\lfloor\frac{#1}{2}\rfloor}}

\def\qpmx#1#2{\begin{pmatrix}#1\\ #2\end{pmatrix}}
\def\qptmx#1#2{\left(\begin{smallmatrix}#1\\ #2\end{smallmatrix}\right)}

\newtheorem{prop}{Proposition}[section]
\newtheorem{proposition}[prop]{Proposition}

\newtheorem{theorem}[prop]{Theorem}
\newtheorem{lemma}[prop]{Lemma}

\newtheorem{corollary}[prop]{Corollary}

\theoremstyle{definition}

\newtheorem{definition}[prop]{Definition}

\newtheorem{notation}[prop]{Notation}

\theoremstyle{remark}
\newtheorem{remark}[prop]{Remark}

\newtheorem*{example*}{Example}

\begin{document}

\begin{abstract}
Trace diagrams are structured graphs with edges labeled by matrices. Each diagram has an interpretation as a particular multilinear function. We provide a rigorous combinatorial definition of these diagrams using a notion of signed graph coloring, and prove that they may be efficiently represented in terms of matrix minors. Using this viewpoint, we provide new proofs of several standard determinant formulas and a new generalization of the Jacobi determinant theorem.
\end{abstract}
%

\maketitle


\section{Introduction}

Trace diagrams provide a graphical means of performing computations in multilinear algebra. The following example, which proves aa vector identity, illustrates the power of the notation.
\begin{example*}
    For $\bfu,\bfv,\bfw\in\C^3$, diagrams for the cross product and inner product are
	$$\bfu\times\bfv =
    \tikz{\node[vertex]at(0,.6){}edge(0,1)edge[bend right]node[vector,pos=1]{$\bfu$}(-.5,0)edge[bend left]node[vector,pos=1]{$\bfv$}(.5,0);}
	\quad\text{and}\quad
	\bfu\cdot\bfv =
    \tikz{\draw(0,0)node[vector]{$\bfu$}to[out=90,in=90,looseness=2](1,0)node[vector]{$\bfv$};}.
    $$
    By ``bending'' the diagrammatic identity
    \begin{equation}\label{eq:3binor}
    \tikz[heightoneonehalf]{\node[vertex]at(0,.5){}edge[bend left](.5,0)edge[bend right](-.5,0);
        \node[vertex]at(0,1){}edge[bend right](.5,1.5)edge[bend left](-.5,1.5)edge(0,.5);}
    =\tikz{\draw(0,0)to[bend left,looseness=.1](.6,1);\draw(.6,0)to[bend right,looseness=.1](0,1);}
    -\tikz{\draw(0,0)to[bend left,looseness=.1](0,1);\draw(.4,0)to[bend right,looseness=.1](.4,1);},
    \end{equation}
    and attaching vectors, one obtains
    $$
    \tikz[scale=.6]{
        \node[vertex](node)at(.5,1){}
            edge[bend right]node[vector,pos=1]{$\bf u$}(0,0)
            edge[bend left]node[vector,pos=1]{$\bf v$}(1,0);
        \node[vertex](node2)at(2.5,1){}
            edge[bend right]node[vector,pos=1]{$\bf w$}(2,0)
            edge[bend left]node[vector,pos=1]{$\bf x$}(3,0)
            edge[bend right=90](node);
    }
    =
    \tikz[scale=.6]{
        \draw(0,0)node[vector]{$\bf u$}to[bend left=90,looseness=2](2,0)node[vector]{$\bf w$};
        \draw(1,0)node[vector]{$\bf v$}to[bend left=90,looseness=2](3,0)node[vector]{$\bf x$};
    }
    -\tikz[scale=.6]{
        \draw(0,0)node[vector]{$\bf u$}to[bend left=90,looseness=1.7](3,0)node[vector]{$\bf x$};
        \draw(1,0)node[vector]{$\bf v$}to[bend left=90,looseness=2.5](2,0)node[vector]{$\bf w$};
    },
    $$
    which is the vector identity
	$$(\bfu\times\bfv)\cdot(\bfw\times\mathbf{x})
        =(\bfu\cdot\bfw)(\bfv\cdot\mathbf{x})-(\bfu\cdot\mathbf{x})(\bfv\cdot\bfw).$$
    We will later prove \eqref{eq:3binor} and show that every step here can be mathematically rigorous. 
\end{example*}

In this paper, we define a set of combinatorial objects called \emph{trace diagrams}. Each diagram translates to a well-defined multilinear function, provided it is \emph{framed} (the framing specifies the domain and range of the function). We introduce the idea of \emph{signed graph coloring} to describe this translation, and show that it preserves a tensorial structure. We prove two results regarding the relationship between multilinear algebra and trace diagrams. Under traditional notation, a multilinear function is characterized by its action on a basis of tensor products in the domain. Theorem \ref{thm:relations} shows that trace diagram notation is more powerful than this standard notation for functions, since a single diagrammatic identity may simultaneously represent several different identities of multilinear functions. In the above example, the diagrammatic identity \eqref{eq:3binor} is used to prove a vector identity; another vector identity arising from the same diagram is given in section \ref{s:functor}.

Our main results concern the ``structural'' properties of trace diagrams. In particular, we characterize their decomposition into \emph{diagram minors}, which are closely related to matrix minors. Theorem \ref{thm:diagram-minors} describes the condition under which this decomposition is possible, and Theorem \ref{thm:diagram-minors-2} gives an upper bound for the number of matrix minors required in a formula for a trace diagram's function.

As an application, we use trace diagrams to provide new proofs of classical determinant identities. Cayley, Jacobi, and other 19th century mathematicians described several methods for calculating determinants in general and for special classes of matrices \cite{Mu1882}. The calculations could often take pages to complete because of the complex notation and the need to keep track of indices. In contrast, we show that diagrammatic proofs of certain classic results come very quickly, once the theory has been suitably developed. One can easily generalize the diagrammatic identities by adding additional matrices, which is not as easy to do with the classical notation for matrices. Theorem \ref{thm:generalized-jacobi}, a novel generalization of a determinant theorem of Jacobi, is proven in this manner.

While the term \emph{trace diagrams} is new, the idea of using diagrammatic notations for algebraic calculations has a rich history \cite{Ba96,Bu97,Cv08,LP09,Ste90}. In the early 1950s, Roger Penrose invented a diagrammatic notation that streamlined calculations in multilinear algebra. In his context, indices became labels on edges between ``spider-like'' nodes, and tensor contraction meant gluing two edges together \cite{Pen71}. In knot theory, Kauffman generalized Penrose's diagrams and described their relation to knot polynomials \cite{Kau91}. Przytycki and others placed Kauffman's work in the context of skein modules \cite{BFK96,Prz91}. The concept of \emph{planar algebras} \cite{Jon99} unifies many of the concepts underlying diagrammatic manipulations. More recently, Kuperberg introduced spiders \cite{Kup96} as a means of studying representation theory. In mathematical physics, Levinson pioneered the use of diagrams to study angular momentum \cite{Lev56}. This approach proved to be extremely useful, with several textbooks written on the topic. Work on these notations and their broader impact on fundamental concepts in physics culminated in books by Stedman \cite{Ste90} and Cvitanovi{\'c} \cite{Cv08}.

The name ``trace diagrams'' was first used in \cite{Pet06} and \cite{LP09}, where diagrams were used to write down an additive basis for a certain ring of invariants. Special cases of trace diagrams have appeared before in the above works, but they are generally used only as a tool for algebraic calculation. This paper differs in emphasizing the diagrams themselves, their combinatorial construction, and their structural properties.

\subsubsection*{This paper is organized as follows.} Section \ref{s:multilinear} provides a short review of multilinear algebra. In Section \ref{s:colorings} we introduce the idea of \emph{signed graph coloring}, which forms the basis for the translation between trace diagrams and multilinear algebra described rigorously in sections \ref{s:trace-diagrams} and \ref{s:functor}. Section \ref{s:building-blocks} describes the basic properties of trace diagram, and section \ref{s:minors} focuses on the fundamental relationship between matrix minors and trace diagram functions. New proofs of classical determinant results are derived in section \ref{s:proofs}. Finally, in section \ref{s:generalization} we prove a new multilinear algebra identity using trace diagrams.

\ifnum0=1
There are many algorithms for computing the determinant of a matrix: directly from the definition, cofactor expansion, Laplace expansion, Dodgson condensation, and more. The intent of these rules is to simplify "by-hand" calculations of the determinant by exhibiting an easy-to-follow algorithm that does not simply rely on writing down $n!$ summands.

Although these techniques are hundreds of years old, they can all be expressed and proven quite simply using a little-known diagrammatic notation. The diagrammatic notation provides a direct means of comparing the various algorithms, and can also provide ways to construct new algorithms for computing the determinant.

This paper explores the connections between the somewhat diverse areas of trace diagrams (a spinoff of "spin networks"), determinants, and "signed" chromatic polynomials. We provide the specific link between these areas and demonstrate how signed chromatic polynomials are quite useful in calculations of trace diagrams; more than that, the entire theory of trace diagrams can be expressed in terms of such combinatorial arguments.

The central idea is that any trace diagram can be computed by simply enumerating all possible colorings of the diagram. To review, an $n$-trace diagram is a graph with vertices of degree 1 or n only. The diagram may be understood as a multilinear function. Under this interpretation, the edges represent vectors in an underlying space, while the vertices in such diagrams encode an anti-symmetry among the attached edges. The vertices of degree 1 represent the inputs and outputs of the function.

With slightly more structure added to the graph, the diagrams became an equivalent category to multilinear algebra, and some calculations are much more easily done using these doodles.

Combinatorics comes into play when we start looking at all possible ways to color the edges of such graphs with $n$ colors, such that no two edges adjacent to the same vertex have the same color. This is somewhat equivalent to an anti-symmetry condition, which causes a diagram with two identical inputs to a vertex to evaluate to zero. In this way, each diagram has a \emph{coloring space} generated by all possible colorings.

What is most surprising is the sheer power afforded by this diagrammatic notation. They turn out to be a very natural language for dealing with concepts in linear algebra such as the determinant, the adjugate, matrix minors, and the characteristic polynomial. They also reduce classical results such as Jacobi's Theorem and its consequence in Dodgson Condensation to trivial one-line proofs. In contrast, most traditional proofs require several pages.
\fi

\section{Multilinear Algebra}\label{s:multilinear}
This section reviews multilinear algebra and tensors. For further reference, a nice introductory treatment of tensors is given in Appendix B of \cite{FH91}.

Let $V$ be a finite-dimensional vector space over a field $\F$. Informally, a \emph{2-tensor} consists of finite sums of vector pairs $(\bfu,\bfv)\in V\times V$ modulo the relations
	$$(\lambda\bfu,\bfv)=\lambda(\bfu,\bfv)=(\bfu,\lambda\bfv)$$
for all $\lambda\in \F$. The resulting term is denoted $\bfu\otimes\bfv$. More generally, a \emph{$k$-tensor} is an equivalence class of $k$-tuples of vectors, where $k$-tuples are equivalent if and only if they differ by the positioning of scalar constants. In other words, if $\prod_{i=1}^k \lambda_i = \prod_{i=1}^k \mu_i = \Lambda$ then
    $$\lambda_1\bfu_1\otimes\cdots\otimes\lambda_k\bfu_k = \mu_1\bfu_1\otimes\cdots\otimes\mu_k\bfu_k = \Lambda \left(\bfu_1\otimes\cdots\otimes\bfu_k\right).$$

Let $N=\{1,2,\ldots,n\}$. In what follows, we assume that $V$ has basis $\{\bse_1,\bse_2,\ldots,\bse_n\}$. The space of $k$-tensors $V^{\otimes k}\equiv V\otimes\cdots\otimes V$ is itself a vector space with $n^k$ basis elements of the form
    $$\bse_{\alpha} \equiv \bse_{\alpha_1} \otimes \bse_{\alpha_2} \otimes \cdots \otimes \bse_{\alpha_k},$$
one for each $\alpha=(\alpha_1,\alpha_2,\ldots,\alpha_k) \in N^k$. By convention $V^{\tensor 0}=\F$.

Let $\langle\cdot,\cdot\rangle$ be the inner product on $V$ defined by $\langle\bse_i,\bse_j\rangle=\delta_{ij},$ where $\delta_{ij}$ is the Kronecker delta. This extends to an inner product on $V^{\otimes k}$ with
    $$\langle\bse_{\alpha},\bse_{\beta}\rangle
        =\delta_{\alpha_1\beta_1}\delta_{\alpha_2\beta_2}\cdots\delta_{\alpha_k\beta_k},$$
making $\{\bse_{\alpha} : \alpha \in N^k\}$ an orthonormal basis for $V^{\otimes k}$.

Given another vector space $W$ over $\F$, a \emph{multilinear function} $f:V^{\otimes k}\to W$ is one that is linear in each term, so that
    $$f\left((\lambda\bfu+\mu\bfv)\otimes\bfu_2\otimes\cdots\otimes\bfu_k\right)
    =\lambda f(\bfu\otimes\bfu_2\otimes\cdots\otimes\bfu_k)
    + \mu f(\bfv\otimes\bfu_2\otimes\cdots\otimes\bfu_k),$$
and a similar identity holds for each of the other $(k-1)$ terms.

Denote by $\mathrm{Fun}(V^{\otimes j},V^{\otimes k})$ the space of \emph{multilinear functions} from $V^{\otimes j}$ to $V^{\otimes k}$. There are two standard ways to combine these functions. First, given $f\in\mathrm{Fun}(V^{\otimes j},V^{\otimes k})$ and $g\in\mathrm{Fun}(V^{\otimes k},V^{\otimes m})$, one may define a composition $g\circ f$. Second, given $f_1\in\mathrm{Fun}(V^{\otimes j_1},V^{\otimes k_1})$ and $f_2\in\mathrm{Fun}(V^{\otimes j_2},V^{\otimes k_2})$, then $f_1\otimes f_2\in\mathrm{Fun}(V^{\otimes(j_1+j_2)},V^{\otimes(k_1+k_2)})$ is the multilinear function defined by letting $f_1$ operate on the first $j_1$ tensor components of $V^{\otimes(j_1+j_2)}$ and $f_2$ on the last $j_2$ components.

A multilinear function $f\in\mathrm{Fun}(V^{\otimes k})\equiv\mathrm{Fun}(V^{\otimes k},\F)$ is commonly called a \emph{multilinear form}. Also, functions $f:\F\to \F$ may be thought of as elements of $\F$. In particular, $\mathrm{Fun}(\F,\F)\cong \F$ via the isomorphism $f\mapsto f(1)$.

The space of tensors $V^{\otimes k}$ is isomorphic to the space of forms $\mathrm{Fun}(V^{\otimes k})$. Given $f\in\mathrm{Fun}(V^{\otimes k})$, the isomorphism maps
    \begin{equation}\label{eq:form-to-tensor}
        f \mapsto \sum_{\alpha\in N^k} f(\bse_{\alpha}) \bse_{\alpha} \in V^{\otimes k}.
    \end{equation}
This is the \emph{duality} property of tensor algebra. Loosely speaking, multilinear functions do not distinguish between inputs and outputs; up to isomorphism all that matters is the total number of inputs and outputs.

One relevant example is the determinant, which can be written as a multilinear function $V^{\otimes k}\to \F$. In particular, if a $k\times k$ matrix is written in terms of its column vectors as $A=\left[\bfa_1 \: \bfa_2 \: \cdots \: \bfa_k\right]$, then the determinant maps the ordered $k$-tuple $\bfa_1\otimes\cdots\otimes\bfa_k$ to $\det(A)$. This may be defined on the tensor product since a scalar multiplied on a single column may be factored outside the determinant. Determinants additionally are \emph{antisymmetric}, since switching any two columns changes the sign of the determinant. Antisymmetric functions can also be considered as functions on an \emph{exterior (wedge) product} of vector spaces, which we do not define here.

\section{Signed Graph Coloring}\label{s:colorings}
This section introduces graph theoretic principles that will be used in defining trace diagram functions. Although the terminology of colorings is borrowed from graph theory, to our knowledge the notion of signed graph coloring is new, being first described in \cite{Pet06}. Some readers may wish to consult a graph theory text such as \cite{Wes01} for further background on graph theory and edge-colorings, or an abstract algebra text such as \cite{Fra02} for further background on permutations.

\subsection{Ciliated Graphs and Edge-Colorings}

A \emph{graph} $G=(V,E)$ consists of a finite collection of \emph{vertices} $V$ and a finite collection of \emph{edges} $E$. Throughout this paper, we permit an edge to be any one of the following:
\begin{enumerate}
    \item a 2-vertex set $\{v_1,v_2\}\subset V$, representing an (undirected) edge connecting vertices $v_1$ and $v_2$;
    \item a 1-vertex set $\{v\}\subset V$ called a \emph{loop}, representing an edge connecting a vertex to itself; or
    \item the empty set $\{\}\subset V$, denoted $\bigcirc$, representing a trivial loop that does not connect any vertices.
\end{enumerate}
In addition, we allow the collection of edges $E$ to contain repeated elements of the same form.

Two vertices are \emph{adjacent} if there is an edge connecting them; two edges are \emph{adjacent} if they share a common vertex. An edge is \emph{adjacent} to a vertex if it contains that vertex. Given a vertex $v$, the set of edges adjacent to $v$ will be denoted $E(v)$. The \emph{degree} $\deg(v)$ of a vertex $v$ is the number of adjacent edges, where any loops at the vertex are counted twice. Vertices of degree 1 are commonly called \emph{leaves}.

\begin{figure}[htb]
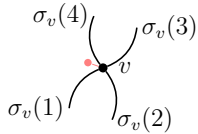

    $$\tikz[shift={(0,.4)}]{
        \node[ciliatednode=170](va)at(0,0){}
            edge[bend right](50:1)edge[bend left](100:1)edge[bend right](230:1)edge[bend left](280:1);
        \node[rightlabel]at(50:1){$\sigma_v(3)$};\node[leftlabel]at(100:1){$\sigma_v(4)$};
        \node[leftlabel]at(230:1){$\sigma_v(1)$};\node[rightlabel]at(280:1){$\sigma_v(2)$};
        \node[rightlabel]at(.2,0){$v$};
    }$$
    \caption{Proceeding counterclockwise from the ciliation at the vertex $v$, one obtains the ordering of edges $\sigma_v(1),\sigma_v(2),\sigma_v(3),\sigma_v(4)$.}\label{f:ciliation}
\end{figure}

\begin{definition}
    A \emph{ciliated graph} $G=(V,E,\sigma_*)$ is a graph $(V,E)$ together with an ordering $\sigma_v:\{1,2,\ldots,\deg(v)\}\to E(v)$ of edges at each vertex $v\in V$.
\end{definition}
By convention, when such graphs are drawn in the plane, the ordering is specified by enumerating edges in a counter-clockwise fashion from a ciliation, as shown in Figure \ref{f:ciliation}.

\begin{definition}\label{d:edge-coloring}
    Given the set $N=\{1,2,\ldots,n\}$, an \emph{$n$-edge-coloring} of a graph $G=(V,E)$ is a map $\kappa:E\to N$. The coloring is said to be \emph{proper} if the graph does not contain any loops and no two adjacent edges have the same label; equivalently, for every vertex $v$ the restriction $\kappa:E(v)\to N$ is one-to-one. When $n$ is clear from context, we denote the set of all proper $n$-edge-colorings of a graph $G$ by $\mathsf{col}(G)$.
\end{definition}
Note that some graphs do not have proper $n$-edge-colorings for certain $n$. As a simple example, the graph
    $\tikz[scale=.8,shift={(0,.2)}]{\draw(0,0)node[vertex]{}to(1,0)node[vertex]{}
      to(.5,.866)node[vertex]{}to(0,0);}$
has no 2-edge-colorings.

\subsection{Permutations and Signatures of Edge-Colorings}

Let $S_n$ denote the set of permutations of $N=\{1,2,\ldots,n\}$. We denote a specific permutation as follows: $\qptmx{1&2&3}{1&2&3}$ denotes the identity permutation, and $\qptmx{1&2&3}{3&2&1}$ denotes the permutation mapping $1\mapsto 3$, $2\mapsto 2$, and $3\mapsto 1$. The \emph{signature} of a permutation is $(-1)^k$, where $k$ is the number of transpositions (or swaps) that must be made to return the permutation to the identity. For example, the permutation $\tbinom{1\:2\:3\:4}{2\:4\:1\:3}$ has signature $-1$, since it takes 3 transpositions to return it to the identity: $$(2,4,1,3)\rightsquigarrow(1,4,2,3)\rightsquigarrow(1,2,4,3)\rightsquigarrow(1,2,3,4).$$

Proper edge-colorings induce permutations at the vertices of ciliated graphs. Given a proper $n$-edge-coloring $\kappa$ and a degree-$n$ vertex $v$, there is a well-defined permutation $\pi_\kappa(v)\in S_n$ defined by
    $$\pi_\kappa(v): i \mapsto \kappa(\sigma_v(i)).$$
In other words, 1 is taken to the label on the first edge adjacent to the vertex, 2 is taken to the label on the second edge, and so on. An example is shown in Figure \ref{f:ciliation2} below.

\begin{definition}[\cite{Pet06}]\label{d:signature-coloring}
    Given a proper $n$-edge-coloring $\kappa$ of a ciliated graph $G=(V,E,\sigma_*)$, the \emph{signature} $\sgn_\kappa(G)$ is the product of permutation signatures on the degree-$n$ vertices:
        $$\sgn_\kappa(G) = \prod_{v\in V_n} \sgn(\pi_\kappa(v)),$$
    where $V_n$ is the set of degree-$n$ vertices in $V$ and $\sgn(\pi_\kappa(v))$ is the signature of the permutation $\pi_\kappa(v)$.
    If there are no degree-$n$ vertices, the signature is $+1$.

    The \emph{signed chromatic index} $\chi(G)$ is the sum of signatures over all proper edge-colorings:
    $$\chi(G) = \sum_{\kappa\in \mathsf{col}(G)} \sgn_\kappa(G).$$
\end{definition}

\begin{figure}[htb]
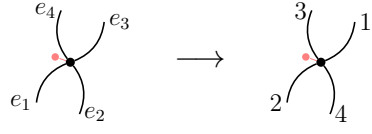

    $$\tikz[shift={(0,.4)}]{
        \node[ciliatednode=170](va)at(0,0){}
            edge[bend right](50:1)edge[bend left](100:1)edge[bend right](230:1)edge[bend left](280:1);
        \node[rightlabel]at(50:1){$e_3$};\node[leftlabel]at(100:1){$e_4$};
        \node[leftlabel]at(230:1){$e_1$};\node[rightlabel]at(280:1){$e_2$};
    }
    \quad\longrightarrow\quad
    \tikz[shift={(0,.4)}]{
        \node[ciliatednode=170](va)at(0,0){}
            edge[bend right](50:1)edge[bend left](100:1)edge[bend right](230:1)edge[bend left](280:1);
        \node[rightlabel]at(50:1){$1$};\node[leftlabel]at(100:1){$3$};
        \node[leftlabel]at(230:1){$2$};\node[rightlabel]at(280:1){$4$};
    }$$
    \caption{The proper edge-coloring at right induces the permutation $\tbinom{1\:2\:3\:4}{2\:4\:1\:3}$ on the ciliated vertex shown. The signature of the coloring is $-1$.}\label{f:ciliation2}
\end{figure}

\begin{example*}
    For $n=2$, the ciliated graph
        $G=
		\tikz[scale=.75,shift={(0,.125)}]{
            \node[ciliatednode=185](va)at(0,0){};
            \node[ciliatednode=175](vb)at(0,1){}edge[bend left=60](va)edge[bend right=60](va);
            \node[rightlabel]at(0,-.2){$w$};
            \node[rightlabel]at(0,1.2){$v$};
        }$
    has exactly two proper edge-colorings:
    \begin{equation}\label{eq:coloring-ex}
        \kappa_1\leftrightarrow
		\tikz{
            \node[ciliatednode=185](va)at(0,0){};
            \node[rightlabel]at(0,-.2){$w$};
            \node[rightlabel]at(0,1.2){$v$};
            \node[ciliatednode=175](vb)at(0,1){}
                edge[bend left=60]node[rightlabel,pos=.5]{$1$}(va)
                edge[bend right=60]node[leftlabel,pos=.5]{$2$}(va);
        }
        \qquad\text{and}\qquad
        \kappa_2\leftrightarrow
		\tikz{
            \node[ciliatednode=185](va)at(0,0){};
            \node[rightlabel]at(0,-.2){$w$};
            \node[rightlabel]at(0,1.2){$v$};
            \node[ciliatednode=175](vb)at(0,1){}
                edge[bend left=60]node[rightlabel,pos=.5]{$2$}(va)
                edge[bend right=60]node[leftlabel,pos=.5]{$1$}(va);
        }.
    \end{equation}
    With the counter-clockwise ordering, $\pi_{\kappa_1}(w)=\qptmx{1&2}{1&2}$ and $\pi_{\kappa_1}(v)=\qptmx{1&2}{2&1}$, so the signature of the first coloring is
        $$\sgn_{\kappa_1}(G)=\sgn(\pi_{\kappa_1}(w))\:\sgn(\pi_{\kappa_1}(v))=\sgn\qptmx{1&2}{1&2}\:\sgn\qptmx{1&2}{2&1}=-1.$$
    In the second case, the permutations are $\qptmx{1&2}{2&1}$ at $w$ and $\qptmx{1&2}{1&2}$ at $v$, so the signature is again $-1$. Therefore, the signed chromatic index of this ciliated graph is $\chi(G)=-2$.
\end{example*}

\subsection{Pre-Edge-Colorings}

\begin{definition}[\cite{Pet06}]\label{d:pre-edge-coloring}
    A \emph{pre-edge-coloring} of a graph $G=(E,V)$ is an edge-coloring $\check\kappa:\check E\to N$ of a subset $\check E\subset E$ of the edges of $G$. A \emph{leaf-coloring} is a pre-edge-coloring of the edges adjacent to the degree-1 vertices.

    Two pre-edge-colorings $\check\kappa_1:\check E_1\to N$ and $\check\kappa_2:\check E_2\to N$ are \emph{compatible} if they agree on the intersection $\check E_1\cap\check E_2$. In this case, the map $\check\kappa_1\cup\check\kappa_2$ defined by $(\check\kappa_1\cup\check\kappa_2)|_{\check E_i} = \check\kappa_i|_{\check E_i}$ is also a pre-edge-coloring.

    If $\check\kappa_1:\check E_1\to N$ and $\check\kappa_2:\check E_2\to N$ are compatible and $\check E_1\subset\check E_2$, we say that $\check\kappa_2$ \emph{extends} $\check\kappa_1$ and write $\check\kappa_2\succ\check\kappa_1$. We denote the (possibly empty) set of proper edge-colorings that extend $\check\kappa$ by
        $$\mathsf{col}_{\check\kappa}(G) \equiv \{\kappa\in\mathsf{col}(G) : \kappa\succ\check\kappa\}.$$
    The \emph{signed chromatic subindex} of a pre-edge-coloring $\check\kappa$ is the sum of signatures of its proper extensions:
        $$\chi_{\check\kappa}(G) = \sum_{\kappa\succ\check\kappa} \sgn_\kappa(G).$$
\end{definition}

\begin{example*}
    For $n=3$, the pre-edge-coloring
        $\check\kappa\leftrightarrow
		\tikz[scale=.75,shift={(0,.25)}]{
            \node[ciliatednode=185](va)at(0,.5){}edge[bend left](-.5,1)edge[bend right](.5,1);
            \node[ciliatednode=175](vb)at(0,0){}edge(va)
                edge[bend right]node[leftlabel,pos=1]{$1$}(-.5,-.5)
                edge[bend left]node[rightlabel,pos=1]{$2$}(.5,-.5);
		}$
    extends to exactly two proper edge-colorings:
    \begin{equation}\label{eq:coloring-ex2}
        \kappa_1\leftrightarrow
		\tikz{
            \node[ciliatednode=185](va)at(0,.5){}
                edge[bend left]node[leftlabel,pos=1]{$1$}(-.5,1)
                edge[bend right]node[rightlabel,pos=1]{$2$}(.5,1);
            \node[ciliatednode=175](vb)at(0,0){}edge node[rightlabel]{$3$}(va)
                edge[bend right]node[leftlabel,pos=1]{$1$}(-.5,-.5)
                edge[bend left]node[rightlabel,pos=1]{$2$}(.5,-.5);
		}
        \qquad\text{and}\qquad
        \kappa_2\leftrightarrow
        \tikz{
            \node[ciliatednode=185](va)at(0,.5){}
                edge[bend left]node[leftlabel,pos=1]{$2$}(-.5,1)
                edge[bend right]node[rightlabel,pos=1]{$1$}(.5,1);
            \node[ciliatednode=175](vb)at(0,0){}edge node[rightlabel]{$3$}(va)
                edge[bend right]node[leftlabel,pos=1]{$1$}(-.5,-.5)
                edge[bend left]node[rightlabel,pos=1]{$2$}(.5,-.5);
		}.
    \end{equation}
    One computes the signed chromatic subindex by summing over the signature of each coloring. In the first case,
        $$\sgn_{\kappa_1}(G)=\sgn\qptmx{1&2&3}{1&2&3}\:\sgn\qptmx{1&2&3}{3&2&1}=-1,$$
    where the permutations are read in counter-clockwise order from the vertex. In the second case, the permutations are $\qptmx{1&2&3}{1&2&3}$ and $\qptmx{1&2&3}{3&1&2}$, so
        $$\sgn_{\kappa_2}(G)=\sgn\qptmx{1&2&3}{1&2&3}\:\sgn\qptmx{1&2&3}{3&1&2}=+1.$$
    Summing the two signatures, the signed chromatic subindex is 
        $$\chi_{\check\kappa}(G)=\sgn_{\kappa_1}(G)+\sgn_{\kappa_2}(G)=-1+1=0.$$
\end{example*}

\section{Trace Diagrams}\label{s:trace-diagrams}

Penrose was probably the first to describe how tensor algebra may be performed diagrammatically \cite{Pen71}. In his framework, edges in a graph represent elements of a vector space, and nodes represent multilinear functions. Trace diagrams are a generalization of Penrose's \emph{tensor diagrams}, in which edges may be labeled by matrices and nodes represent the determinant.

The closest concept in traditional graph theory is a \emph{voltage graph} (also called a \emph{gain graph}), in which the edges of a graph are marked by group elements in an ``orientable'' way \cite{Gro74}. Diagrams labeled by matrices also make frequent appearances in skein theory \cite{Bu97,Sik01} and occasional appearance in work of Stedman and Cvitanovi{\'c} \cite{Cv08,Ste90}.

\subsection{Definition}

In the remainder of this paper, $V$ will represent an $n$-dimensional vector space over a base field $\mathbb{F}$ (with $n\ge2$), and $\{\bse_1,\bse_2,\ldots,\bse_n\}$ will represent an orthonormal basis for $V$.

\begin{figure}
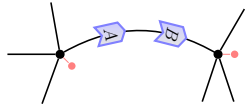

    $$\tikz{
        \node[ciliatednode=-45]at(0,0){}
            edge(100:1)edge(180:1)edge(250:1)
            edge[bend left]node[matrix on edge,pos=.3]{$A$}node[matrix on edge,pos=.7]{$B$}(3,0);
        \node[ciliatednode=0]at(3,0){}edge(3.5,.86)edge(2.7,-.9)edge(3.3,-.9);
    }$$
\caption{An unframed 4-trace diagram. Degree-$n$ vertices are ciliated and degree-2 vertices are marked by matrices in an oriented manner.}\label{f:tracediagram1}
\end{figure}

\begin{definition}\label{d:tracediagrams}
    An \emph{$n$-trace diagram} is a ciliated graph $\mathcal{D}=(V_1\sqcup V_2\sqcup V_n,E,\sigma_*)$, where $V_i$ is comprised of vertices of degree $i$, together with a labeling $A_{\mathcal{D}}:V_2\to\mathrm{Fun}(V,V)$ of degree-2 vertices by linear transformations. If there are no degree-1 vertices, the diagram is said to be \emph{closed}.

    A \emph{framed trace diagram} is a diagram together with a partition of the degree-1 vertices $V_1$ into two disjoint ordered collections: the \emph{inputs} $V_I$ and the \emph{outputs} $V_O$.
\end{definition}
Thus, trace diagrams contain vertices of degree 1, 2, or $n$ only, and the degree-2 vertices represent matrices. An example is shown in Figure \ref{f:tracediagram1}.  Note that in the case $n=2$, the vertices in $V_2$ and $V_n$ have the same degree but are disjoint sets. By convention, framed trace diagrams are drawn with inputs at the bottom of the diagram and outputs at the top. Both are assumed to be ordered left to right.

As shown in Figure \ref{f:tracediagram1}, we represent matrix markings at the degree-2 vertices as follows:
    $$
    A \leftrightarrow \tikz{\draw(0,0)to node[matrix on edge]{$A$}(0,1);},
    \quad\quad
    A^{-1} \leftrightarrow \tikz{\draw(0,0)to node[matrix on edge]{$\overline{A}$}(0,1);}.
    $$
Note that when drawing the inverse of a matrix in a diagram, we use the shorthand $\overline{A}$ because the traditional notation $A^{-1}$ is overly cumbersome.

The ordering at a degree-2 vertex $v$ given by the ciliation is implicit in the orientation of the node. Precisely, the ciliation $\sigma:\{1,2\}\to E(v)$ orders the adjacent edges as follows:
    $$
    \tikz{\draw(0,0)node[rightlabel]{$\sigma_v(1)$}to node[matrix on edge]{$A$}(0,1)node[rightlabel]{$\sigma_v(2)$};}.
    $$
We refer the first edge $\sigma_v(1)$ as the ``incoming'' edge and the second edge $\sigma_v(2)$ as the ``outgoing'' edge. In general,
    $\tikz{\draw(0,0)to node[matrix on edge]{$A$}(0,1);}
    \neq\tikz{\draw(0,1)to node[matrix on edge]{$A$}(0,0);}$
since the nodes occur with opposite orientations.

\subsection{Trace Diagram Colorings and their Coefficients}

Trace diagrams require a slightly different definition of edge-coloring:
\begin{definition}\label{d:tracediagram-coloring}
    A \emph{coloring} of an $n$-trace diagram $\mathcal{D}$ is a map $\kappa:E\to N$. The coloring is \emph{proper} if the labels at each $n$-vertex are distinct. The (possibly empty) space of all colorings of $\mathcal{D}$ is denoted $\mathsf{col}(\mathcal{D})$.
\end{definition}
Note that in a proper coloring of a trace diagram, the edges adjoining a matrix may have the same label.
\begin{definition}\label{d:tracediagram-coefficient}
    Given a coloring $\kappa$ of a trace diagram $\mathcal{D}$ with matrix labeling $A_{\mathcal{D}}:V_2\to\mathrm{Fun}(V,V)$, the \emph{coefficient} $\psi_\kappa(\mathcal{D})$ of the coloring is defined to be
        $$\psi_\kappa(\mathcal{D})\equiv\prod_{v\in V_2}(A_{\mathcal{D}}(v))_{\sigma_v(2)\sigma_v(1)},$$
    where $(A(v))_{\sigma_v(2)\sigma_v(1)}=\langle\bse_{\sigma_v(2)},A(v)\bse_{\sigma_v(1)}\rangle$ represents the matrix entry in the $\sigma_v(2)$\textit{th} row and $\sigma_v(1)$\textit{th} column.
\end{definition}

\begin{example*}
    In the simplest colored diagram with a matrix,
    \begin{equation}\label{eq:diagram-matrix-entry}
         \psi\left(
         \tikz{\draw(0,-.15)node[rightlabel]{$j$}to node[matrix on edge]{$A$}(0,1.15)node[rightlabel]{$i$};}
         \right) = (A)_{ij}.
    \end{equation}
    Similarly,
    $$
    \psi\left(
    \tikz[heightoneonehalf]{\draw(0,0)to(0,1.5)
        node[rightlabel,pos=0,scale=.8]{$\:\:\:k$}
        node[matrix on edge,scale=.8,pos=.25]{$B$}
        node[rightlabel,pos=.5,scale=.8]{$\:\:\:j$}
        node[matrix on edge,scale=.8,pos=.75]{$A$}
        node[rightlabel,pos=1,scale=.8]{$\:\:\:i$};}
    \right) = (A)_{ij}(B)_{jk}.$$
\end{example*}

\begin{example*}
    In the colored diagram
    $\tikz[heightoneonehalf]{
            \node[ciliatednode=160](va)at(0,1.5){};
            \node[ciliatednode=200](vb)at(0,0){}
                edge[out=25,in=-25]node[matrix on edge]{$A$}
                    node[rightlabel,pos=.95]{$\:1$}node[rightlabel,pos=.05]{$\:2$}(va)
                edge[out=155,in=-155]node[matrix on edge]{$A$}
                    node[leftlabel,pos=.95]{$2\:$}node[leftlabel,pos=.05]{$1\:$}(va);},
    $
    the coefficient is $(A)_{21}(A)_{12}$.
\end{example*}


\subsection{Trace Diagram Functions}

Recall that $\{\bse_1,\ldots,\bse_n\}$ represents an orthonormal basis for the vector space $V$. In a framed trace diagram, a basis element $\bse_{\alpha}\in V^{\tensor|V_I|}$ is equivalent to a labeling of the input vertices by basis elements. This labeling induces a pre-coloring on the adjacent edges: if a vertex is labeled by $\bse_i$, then its adjacent edge is labeled by $i$. We denote this pre-coloring by $\alpha$. Likewise, a basis element $\bse_{\beta}\in V^{\tensor|V_O|}$ induces a pre-coloring on edges adjacent to the output vertices, which we denote $\beta$. Since $V_I$ and $V_O$ comprise all degree-1 vertices in the diagram, if $\alpha$ and $\beta$ exist (and are compatible) then $\alpha\cup\beta$ is a leaf-coloring of the diagram.

The next definition is the key concept relating trace diagrams and multilinear functions. Each diagram corresponds to a unique function, whose coefficients are the signed chromatic subindices of these leaf-colorings, weighted by coloring coefficients.

\begin{definition}\label{d:diagram-weight}
    Given a trace diagram $\mathcal{D}$, the \emph{weight} $\chi_{\gamma}(\mathcal{D})$ of a leaf-coloring $\gamma$ is
    \begin{equation}
        \chi_{\gamma}(\mathcal{D})=\sum_{\kappa\succ\gamma} \sgn_{\kappa}(\mathcal{D})\psi_{\kappa}(\mathcal{D}).
    \end{equation}
    The \emph{value} of a closed diagram $\mathcal{D}$ is
        $$\chi(\mathcal{D})=\sum_{\kappa\in\mathsf{col}(\mathcal{D})}
        \sgn_{\kappa}(\mathcal{D})\psi_{\kappa}(\mathcal{D}).$$
\end{definition}
\begin{definition}\label{d:diagram-function}
    Given a framed trace diagram $\mathcal{D}$, the \emph{trace diagram function} $f_\mathcal{D}: V^{\tensor|V_I|} \to V^{\tensor|V_O|}$ is the linear extension of the basis mappings
    \begin{equation}\label{eq:diagram-function}
        f_{\mathcal{D}}:
        \bse_{\alpha}
        \mapsto \sum_{\beta\in N^{|V_O|}}
            \chi_{\alpha\cup\beta}(\mathcal{D}) \bse_{\beta},
    \end{equation}
    where $f_{\mathcal{D}}:\bse_{\alpha}\mapsto 0$ if $\bse_\alpha$ does not induce a pre-coloring or does not extend to any proper colorings.
\end{definition}

\begin{remark}
    If $n$ is odd, trace diagrams may be drawn without ciliations, since $\sgn(\sigma)$ is invariant under cyclic re-orderings:
    $$\sgn\qpmx{1 & \cdots & n-1 & n}{a_1 & a_2 & \cdots & a_n} = \sgn\qpmx{1 & \cdots & n-1 & n}{a_2 & \cdots & a_n & a_1}.$$
\end{remark}

We will sometimes abuse notation by using the diagram $\mathcal{D}$ interchangeably with $f_\mathcal{D}$.
When describing a diagram's function, we will sometimes mark the input vertices by vectors to indicate the input vectors. For example,
    $$
    
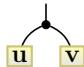
\tikz{\node[vertex]at(0,.6){}edge(0,1)edge[bend right]node[vector,pos=1]{$\bfu$}(-.5,0)edge[bend left]node[vector,pos=1]{$\bfv$}(.5,0);}
    $$
is used as shorthand for $f_{\tikz[scale=.5,shift={(0,.2)}]{\node[vertex]at(0,.5){}edge(0,1)edge[bend right](-.5,0)edge[bend left](.5,0);}}(\bfu\otimes\bfv)$.
We also write formal linear sums of diagrams to indicate the corresponding sums of functions. See the next section for explicit details on why this is permissible.

\subsection{Computations and Examples}

The next few examples show how to compute the value of a closed diagram. Later examples will demonstrate how trace diagram functions are computed.

\begin{example*}
    The ``barbell'' diagram
    \tikz[rotate=90,shift={(.5,0)}]{\draw(0,0)circle(.25);\draw(0,.75)circle(.25);\draw(0,.25)to(0,.5);}
    has no proper colorings, since in any coloring the same color meets a vertex twice. Therefore, the diagram's value is
    $\chi(\tikz[rotate=90,shift={(.5,0)}]{\draw(0,0)circle(.25);\draw(0,.75)circle(.25);\draw(0,.25)to(0,.5);})=0$.
\end{example*}

\begin{example*}
    The simple loop $\bigcirc$ (with no vertices) has $n$ proper colorings, $\left\{\bigcirc^i\right\}$ for $i=1,2,\ldots,n$. Since there are no vertices, the weight of each coloring is $+1$. Hence, the value of the circle is $\chi(\bigcirc)=\sum_{i=1}^n 1 = n$.
\end{example*}

The next example is the reason for the terminology `trace' diagram.

\begin{example*}
    The simplest closed trace diagram with a matrix is
        \tikz[scale=.8,shift={(0,.1)}]{\draw(0,.5)circle(.5);\node[small matrix,yscale=.8]at(.5,.5){$A$};}.
    There are $n$ proper colorings of the form
        \tikz[scale=.8,shift={(0,.1)}]{\draw(0,.5)circle(.5);\node[small matrix,yscale=.8]at(.5,.5){$A$};
            \node[basiclabel]at(.5,-.1){$i$};},
    for $i=1,2,\ldots,n$. The coefficient of the $i$th coloring is $(A)_{ii}\equiv a_{ii}$, so the diagram's value is
    \begin{equation}\label{eq:trace-diagram}
        \chi\left(\tikz{\draw(0,.5)circle(.5);\node[small matrix,yscale=.8]at(.5,.5){$A$};}\right)
        =a_{11}+\cdots+a_{nn}=\tr(A).
    \end{equation}
\end{example*}

The propositions that follow will be used later in this paper, but they are also intended as examples illustrating how to compute trace diagram functions.

\begin{proposition}\label{p:diagram-identity}
    The function of the diagram $\tikz{\draw(0,0)to[wavyup](.05,1);}$ is the identity $\bfv\mapsto\bfv$.
\begin{proof}
    To compute $f_{|}(\bse_i)$, one considers the pre-coloring $\alpha$ in which the input edge has been labeled $i$.  But this is also a full coloring, and since there are no vertices and no matrices, the weight of that coloring is $+1$. Hence, $\beta=\alpha=(i)$ is the only summand in \eqref{eq:diagram-function} and $f_{|}(\bse_i)=\bse_i$. By linear extension, this means $f_{|}(\bfv)=\bfv$ for all $\bfv\in V$, so the diagram's function is the identity on $V$.
\end{proof}
\end{proposition}

\begin{proposition}\label{p:mxmult}
    Given $n\times n$ matrices $A$ and $B$, (i)
        $f_{\tikz[scale=.65]{\draw(0,0)to[wavyup](.05,1)node[matrix on edge,pos=.5,scale=.8]{$A$};}}
            :\bfv\mapsto A\bfv$
    and
    (ii) the diagrams
        \tikz[heightoneonehalf]{
            \draw(0,0)to(0,1.5)
                node[matrix on edge,pos=.25]{$B$}
                node[matrix on edge,pos=.75]{$A$};}
    and
        \tikz{\draw(0,0)to[wavyup](.05,1)node[matrix on edge,pos=.5]{$AB$};}
    have the same function.
\begin{proof}
    Recall that the coefficient of a coloring of
        \tikz{\draw(0,0)to[wavyup](.05,1)node[matrix on edge,pos=.5]{$A$};}
    is $(A)_{ij}$, where $i$ is the label at the top of the diagram and $j$ is the label at the bottom of the diagram \eqref{eq:diagram-matrix-entry}. Thus,
        $$f_{\tikz[scale=.65]{\draw(0,0)to[wavyup](.05,1)node[matrix on edge,pos=.5,scale=.8]{$A$};}}
            :\bse_j \mapsto
              \sum_{i=1,\ldots,n} \psi\left(
                    \tikz{\draw(0,-.15)node[rightlabel]{$j$}to node[matrix on edge]{$A$}(0,1.15)node[rightlabel]{$i$};}
                \right)
            = \sum_{i=1,\ldots,n} (A)_{ij} \equiv A\bse_j.$$
    By linear extension,
        $f_{\tikz[scale=.65]{\draw(0,0)to[wavyup](.05,1)node[matrix on edge,pos=.5,scale=.8]{$A$};}}
            :\bfv \mapsto A\bfv$,
    verifying the first result.

    In the case of the diagram
        $\tikz[heightoneonehalf]{
            \draw(0,0)to(0,1.5)
                node[matrix on edge,pos=.25]{$B$}
                node[matrix on edge,pos=.75]{$A$};}$
    with two matrices, one reasons similarly to show that the diagram's function maps
        $$\bse_k \longmapsto
            \sum_{i=1,\ldots,n}\sum_{j=1,\ldots,n} \psi\left(
                \tikz[heightoneonehalf]{\draw(0,0)to(0,1.5)
                    node[rightlabel,pos=0,scale=.8]{$\:\:\:k$}
                    node[matrix on edge,scale=.8,pos=.25]{$B$}
                    node[rightlabel,pos=.5,scale=.8]{$\:\:\:j$}
                    node[matrix on edge,scale=.8,pos=.75]{$A$}
                    node[rightlabel,pos=1,scale=.8]{$\:\:\:i$};}
            \right) =
            \sum_{i=1,\ldots,n}\sum_{j=1,\ldots,n} (A)_{ij}(B)_{jk} \equiv AB\bse_k.$$
    Thus, $\bfv\mapsto (AB)\bfv$, verifying the second result.
\end{proof}
\end{proposition}

We can now prove the diagrammatic identity \eqref{eq:3binor} stated in the introduction.
\begin{proposition}\label{p:3binor}
    As a statement about the functions underlying the corresponding 3-trace diagrams,
    $$
    \tikz[heightoneonehalf]{\node[vertex]at(0,.5){}edge[bend left](.5,0)edge[bend right](-.5,0);
        \node[vertex]at(0,1){}edge[bend right](.5,1.5)edge[bend left](-.5,1.5)edge(0,.5);}
    =\tikz{\draw(0,0)to[bend left,looseness=.1](.6,1);\draw(.6,0)to[bend right,looseness=.1](0,1);}
    -\tikz{\draw(0,0)to[bend left,looseness=.1](0,1);\draw(.4,0)to[bend right,looseness=.1](.4,1);}.
    $$
\begin{proof}
    Proposition \ref{p:diagram-identity} implies that
    $$\tikz{\draw(0,0)to[bend left,looseness=.1](.6,1);\draw(.6,0)to[bend right,looseness=.1](0,1);}
        :\bfu\otimes\bfv\mapsto\bfv\otimes\bfu
    \qquad\text{and}\qquad
    \tikz{\draw(0,0)to[bend left,looseness=.1](0,1);\draw(.4,0)to[bend right,looseness=.1](.4,1);}
        :\bfu\otimes\bfv\mapsto\bfu\otimes\bfv.$$

    Now consider the function for the 3-diagram
        $\mathcal{D}=\tikz[scale=.65,shift={(0,.5)}]{
            \node[ciliatednode=185](va)at(0,.5){}edge[bend left](-.5,1)edge[bend right](.5,1);
            \node[ciliatednode=175](vb)at(0,0){}edge(va)edge[bend right](-.5,-.5)edge[bend left](.5,-.5);
		}$.
    The basis element $\bse_i\otimes\bse_i$, where $i\in\{1,2,3\}$, corresponds to $\alpha=(i,i)$ and induces the pre-coloring $\alpha\leftrightarrow
		\tikz[scale=.65,shift={(0,.5)}]{
            \node[ciliatednode=185](va)at(0,.5){}edge[bend left](-.5,1)edge[bend right](.5,1);
            \node[ciliatednode=175](vb)at(0,0){}edge(va)
                edge[bend right]node[leftlabel,pos=1]{$i$}(-.5,-.5)
                edge[bend left]node[rightlabel,pos=1]{$i$}(.5,-.5);
		}$,
    which does not extend to any proper colorings. Therefore
        $f_{\mathcal{D}}:\bse_i\tensor\bse_i \mapsto 0.$

    The basis element $\bse_i\otimes\bse_j$, where $i\neq j$, induces the pre-coloring
        $\alpha\leftrightarrow
		\tikz[scale=.65,shift={(0,.5)}]{
            \node[ciliatednode=185](va)at(0,.5){}edge[bend left](-.5,1)edge[bend right](.5,1);
            \node[ciliatednode=175](vb)at(0,0){}edge(va)
                edge[bend right]node[leftlabel,pos=1]{$i$}(-.5,-.5)
                edge[bend left]node[rightlabel,pos=1]{$j$}(.5,-.5);
		}$.
    The summation in \eqref{eq:diagram-function} is nominally over 9 possibilities (the number of elements in $N\times N$), but we only need to consider the two full colorings that extend this pre-coloring. These are
    $$
        \alpha\cup\beta_1\leftrightarrow
		\tikz[shift={(0,.25)}]{
            \node[ciliatednode=185](va)at(0,.5){}
                edge[bend left]node[leftlabel,pos=1]{$i$}(-.5,1)
                edge[bend right]node[rightlabel,pos=1]{$j$}(.5,1);
            \node[ciliatednode=175](vb)at(0,0){}edge node[rightlabel]{$k$}(va)
                edge[bend right]node[leftlabel,pos=1]{$i$}(-.5,-.5)
                edge[bend left]node[rightlabel,pos=1]{$j$}(.5,-.5);
		}
        \qquad\text{and}\qquad
        \alpha\cup\beta_2\leftrightarrow
        \tikz[shift={(0,.25)}]{
            \node[ciliatednode=185](va)at(0,.5){}
                edge[bend left]node[leftlabel,pos=1]{$j$}(-.5,1)
                edge[bend right]node[rightlabel,pos=1]{$i$}(.5,1);
            \node[ciliatednode=175](vb)at(0,0){}edge node[rightlabel]{$k$}(va)
                edge[bend right]node[leftlabel,pos=1]{$i$}(-.5,-.5)
                edge[bend left]node[rightlabel,pos=1]{$j$}(.5,-.5);
		},
    $$
    where $k\in\{1,2,3\}$ is not equal to $i$ or $j$. The signatures are $\sgn_{\alpha\cup\beta_1}(\mathcal{D})=-1$ and $\sgn_{\alpha\cup\beta_2}(\mathcal{D})=+1$. This statement was proven in detail for the case of $i=1$ and $j=2$ in \eqref{eq:coloring-ex2}; the other cases are proven similarly. Since there are no matrices in the diagram, the coefficients of the colorings are both 1, and the weights are equal to the signatures.
    Summing over $\bse_{\beta}$ gives
        $$f_{\mathcal{D}}:\bse_i\tensor\bse_j \mapsto -\bse_i\tensor\bse_j + \bse_j\tensor\bse_i.$$
    Combining this with the fact that $f_{\mathcal{D}}:\bse_i\otimes\bse_i\mapsto 0$ proves the general statement
        $$f_{\mathcal{D}} : \bfu\otimes\bfv \mapsto \bfv\otimes\bfu - \bfu\otimes \bfv,$$
    which completes the proof.
\end{proof}
\end{proposition}

We close this section with the diagrams for the inner and cross products.

\begin{proposition}\label{p:diagram-cap}
    The inner product $\bfu\otimes\bfv \mapsto \bfu\cdot\bfv$ of $n$-dimensional vectors is represented by the $n$-trace diagram
        \tikz{\draw(0,0)to[out=90,in=90,looseness=2](1,0);}.
\begin{proof}
    Since there is only one edge, $\bse_i\otimes\bse_j$ does not induce a coloring unless $i=j$. In this case, the weight of the coloring is 1. Therefore, $\bse_i\otimes\bse_j \mapsto 1$ if $i=j$, or $0$ if $i\neq j$. By extension, $\tikz{\draw(0,0)node[vector]{$\bfu$}to[out=90,in=90,looseness=2](1,0)node[vector]{$\bfv$};}=\bfu\cdot\bfv$.
\end{proof}
\end{proposition}

\begin{proposition}\label{p:diagram-cross}
    The cross product $\bfu\otimes\bfv \mapsto \bfu\times\bfv$ of 3-dimensional vectors is represented by the 3-diagram
        \tikz{\node[vertex]at(0,.6){}edge(0,1)edge[bend right](-.5,0)edge[bend left](.5,0);}.
\begin{proof}
    The input $\bse_i\otimes\bse_j$ corresponds to the pre-coloring
        \tikz{\node[vertex]at(0,.6){}edge(0,1)
            edge[bend right]node[leftlabel,pos=1]{$i$}(-.5,0)
            edge[bend left]node[rightlabel,pos=1]{$j$}(.5,0);}.
    If $i=j$, there is no proper coloring extending this pre-coloring, so the diagram's function maps $\bse_i\otimes\bse_i\mapsto 0$. Otherwise, the only proper coloring is
        \tikz{\node[vertex]at(0,.6){}edge node[rightlabel,pos=1]{$k$}(0,1)
            edge[bend right]node[leftlabel,pos=1]{$i$}(-.5,0)
            edge[bend left]node[rightlabel,pos=1]{$j$}(.5,0);},
    where $k$ is not equal to $i$ or $j$. The signature of this coloring is $\qptmx{1&2&3}{i&j&k}$. Thus, $\bse_i\otimes\bse_j\mapsto\sgn\qptmx{1&2&3}{i&j&k}\bse_k$. It is straightforward to check that this extends to the standard cross product; for instance, $\bse_1\otimes\bse_2\mapsto\sgn\qptmx{1&2&3}{1&2&3}\bse_3=\bse_3$. The other cases are similar.
\end{proof}
\end{proposition}

\subsection{Transpose Diagrams}

Given a trace diagram $\mathcal{D}$, we define the \emph{transpose diagram} $\mathcal{D}^*$ to be the trace diagram in which all orientations of matrix vertices in $\mathcal{D}$ have been reversed. The following result describes the relationship between the functions of $\mathcal{D}$ and $\mathcal{D}^*$.
\begin{proposition}[Transpose Diagrams]\label{p:transpose-diagram}
    Let $\mathcal{D}$ be a trace diagram and let $\mathcal{D}^T$ represent the same diagram in which all matrices have been replaced by their transpose. Then $f_{\mathcal{D}^*}=f_{\mathcal{D}^T}$.
\begin{proof}
    By \eqref{eq:diagram-matrix-entry},
        $$
        \psi\left(
            \tikz{\draw(0,1.15)node[rightlabel]{$i$}to(0,-.15)node[rightlabel]{$j$}node[matrix on edge,pos=.5]{$A$};}
        \right)
        =(A)_{ji}=(A^T)_{ij}
        =\psi\left(
            \tikz{\draw(0,-.15)node[rightlabel]{$j$}to node[matrix on edge]{$A^T$}(0,1.15)node[rightlabel]{$i$};}
        \right).$$
    Thus, the impact of transposing matrices on the underlying function is the same as that of reversing 2-vertex orientations.
\end{proof}
\end{proposition}

\section{Multilinear Functions and Diagrammatic Relations}\label{s:functor}

\subsection{Composition and Tensor Product Diagrams}

\begin{figure}
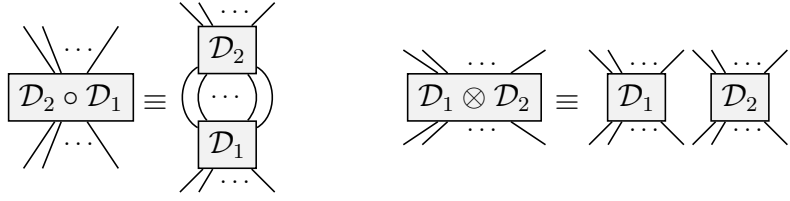

    $$
    \tikz[heighttwo]{
        \node[rectangle,draw=black,fill=black!5]at(0,1){$\mathcal{D}_2\circ\mathcal{D}_1$}
            edge(1,-.5)edge(-.6,-.5)edge(-1,-.5)
            edge(-1,2.5)edge(-.6,2.5)edge(1,2.5);
        \foreach\xa/\xb in{.2/0,.2/2}{\node[basiclabel]at(\xa,\xb){$\cdots$};}
    }
    \equiv
    \tikz[heighttwo]{                
        \node[rectangle,draw=black,fill=black!5](n1)at(0,0){$\mathcal{D}_1$}
            edge(1,-1)edge(-.6,-1)edge(-1,-1);
        \node[rectangle,draw=black,fill=black!5](n2)at(0,2){$\mathcal{D}_2$}
            edge[bend left=40](n1)edge[bend right=40](n1)edge[bend left=60](n1)edge[bend right=60](n1)
            edge(-1,3)edge(-.6,3)edge(1,3);
        \foreach\xa/\xb in{.2/-.8,0/1,.2/2.8}{\node[basiclabel]at(\xa,\xb){$\cdots$};}
    }
    \qquad\qquad
    \tikz[heighttwo]{
        \node[rectangle,draw=black,fill=black!5]at(0,1){$\mathcal{D}_1\otimes\mathcal{D}_2$}
            edge(1.5,0)edge(-1.1,0)edge(-1.5,0)
            edge(-1.5,2)edge(-1.1,2)edge(1.5,2);
        \foreach\xa/\xb in{.2/.3,.2/1.7}{\node[basiclabel]at(\xa,\xb){$\cdots$};}
    }
    \equiv\tikz[heighttwo]{
        \node[rectangle,draw=black,fill=black!5]at(0,1){$\mathcal{D}_1$}
            edge(1,0)edge(-.6,0)edge(-1,0)
            edge(-1,2)edge(-.6,2)edge(1,2);
        \foreach\xa/\xb in{.2/.3,.2/1.7}{\node[basiclabel]at(\xa,\xb){$\cdots$};}
    }\:
    \tikz[heighttwo]{
        \node[rectangle,draw=black,fill=black!5]at(0,1){$\mathcal{D}_2$}
            edge(1,0)edge(-.6,0)edge(-1,0)
            edge(-1,2)edge(-.6,2)edge(1,2);
        \foreach\xa/\xb in{.2/.3,.2/1.7}{\node[basiclabel]at(\xa,\xb){$\cdots$};}
    }
    $$
\caption{The composition of trace diagrams is formed by drawing one diagram above another (left). The tensor product of trace diagrams is found by drawing diagrams side-by-side (right).}\label{f:diagram-operations}
\end{figure}

Given the base field $\F$, let $\mathfrak{D}(I,O)$ denote the free $\F$-module over framed trace diagrams with $I=|V_I|$ inputs and $O=|V_O|$ outputs. There are two ways to combine elements of these spaces. Given $\mathcal{D}_1\in\mathfrak{D}(I_1,O_1)$ and $\mathcal{D}_2\in\mathfrak{D}(I_2,O_2)$ with $|O_1|=|I_2|$, one may form the \emph{composition diagram} $\mathcal{D}_2\circ\mathcal{D}_1$ by gluing the output strands of $\mathcal{D}_1$ to the input strands of $\mathcal{D}_2$. Since by convention inputs are drawn at the bottom of a diagram and outputs at the top, this composition involves drawing one diagram above another. Second, given arbitrary framed diagrams  $\mathcal{D}_1\in\mathfrak{D}(I_1,O_1)$ and $\mathcal{D}_2\in\mathfrak{D}(I_2,O_2)$, we define the \emph{tensor product diagram} $\mathcal{D}_1\otimes\mathcal{D}_2\in\mathfrak{D}(I_1+I_2,O_1+O_2)$ to be that obtained by placing $\mathcal{D}_2$ to the right of $\mathcal{D}_1$. See Figure \ref{f:diagram-operations} for depictions of these two diagram operations.

Both of these structures are preserved under the mapping $\mathcal{D}\mapsto f_{\mathcal{D}}$. The proof is rather technical, but straightforward.

\begin{theorem}\label{thm:categorical}
    Let $\mathcal{D}_1\in\mathfrak{D}(I_1,O_1)$ and $\mathcal{D}_2\in\mathfrak{D}(I_2,O_2)$. The trace diagram function $f_{\mathcal{D}}$ satisfies (i) $f_{\mathcal{D}_1\otimes\mathcal{D}_2}=f_{\mathcal{D}_1}\otimes f_{\mathcal{D}_2}$, and (ii)  $f_{\mathcal{D}_2\circ\mathcal{D}_1}=f_{\mathcal{D}_2}\circ f_{\mathcal{D}_1}$ (when the composition $\mathcal{D}_2\circ\mathcal{D}_1$ is defined).
\begin{proof}
    To see that the tensorial structure is preserved, observe that
        \begin{align*}
        f_{\mathcal{D}_1\otimes\mathcal{D}_2}(\bse_{\alpha_1}\otimes\bse_{\alpha_2})
        &= \sum_{\beta_1,\beta_2\in N^{|O_1|+|O_2|}}
            \chi_{\alpha_1\cup\alpha_2\cup\beta_1\cup\beta_2}(\mathcal{D}_1\otimes\mathcal{D}_2) \bse_{\beta_1}\otimes\bse_{\beta_2}\\
        &= \sum_{\beta_1\in N^{|O_1|}} \sum_{\beta_2\in N^{|O_2|}} \chi_{\alpha_1\cup\beta_1}(\mathcal{D}_1) \chi_{\alpha_2\cup\beta_2}(\mathcal{D}_2)\bse_{\beta_1}\otimes\bse_{\beta_2}\\
        &= \left(\sum_{\beta_1\in N^{|O_1|}} \chi_{\alpha_1\cup\beta_1}(\mathcal{D}_1) \bse_{\beta_1}\right)
        \otimes
            \left(\sum_{\beta_2\in N^{|O_2|}} \chi_{\alpha_2\cup\beta_2}(\mathcal{D}_2) \bse_{\beta_2}\right)\\
        &=f_{\mathcal{D}_1}\otimes f_{\mathcal{D}_2}\:(\bse_{\alpha_1}\otimes\bse_{\alpha_2}).
        \end{align*}

    For composition, assume $\mathcal{D}_2\circ\mathcal{D}_1$ is defined. Apply \eqref{eq:diagram-function} twice to get
        \begin{equation}\label{eq:composition-proof1}
        f_{\mathcal{D}_2}\circ f_{\mathcal{D}_1} : \bse_{\alpha} \mapsto
        \sum_{\gamma\in N^{|O_2|}} \left(\sum_{\beta\in N^{|O_1|}}
            \chi_{\alpha\cup\beta}(\mathcal{D}_1)\chi_{\beta\cup\gamma}(\mathcal{D}_2)\right) \bse_{\gamma}.
        \end{equation}
    The following lemma simplifies the term in parentheses:

    \begin{lemma}
        \begin{equation}\label{eq:composition-proof2}
        \sum_{\beta\in N^{|O_1|}}
            \chi_{\alpha\cup\beta}(\mathcal{D}_1)\chi_{\beta\cup\gamma}(\mathcal{D}_2)
        =\chi_{\alpha\cup\gamma}(\mathcal{D}_2\circ\mathcal{D}_1).
        \end{equation}
    \begin{proof}
        Recall that by definition $\chi_{\alpha\cup\gamma}(\mathcal{D}_2\circ\mathcal{D}_1)$ is defined as a sum over all proper colorings $\kappa$ of the composition diagram $\mathcal{D}_2\circ\mathcal{D}_1$ that extend the pre-coloring $\alpha\cup\gamma$. A proper coloring $\kappa$ induces proper colorings $\kappa_1$ of $\mathcal{D}_1$ and $\kappa_2$ of $\mathcal{D}_2$ that agree on the common edges. So we may write the righthand side of \eqref{eq:composition-proof2} as
            \begin{align*}
            \chi_{\alpha\cup\gamma}(\mathcal{D}_2\circ\mathcal{D}_1)
            &=\sum_{\kappa\succ\alpha\cup\gamma}
                \sgn_\kappa(\mathcal{D}_2\circ\mathcal{D}_1)\psi_\kappa(\mathcal{D}_2\circ\mathcal{D}_1)\\
            &=\sum_{\beta\in N^{|O_1|}} \sum_{\kappa\succ\alpha\cup\beta\cup\gamma} \sgn_\kappa(\mathcal{D}_2\circ\mathcal{D}_1)\psi_\kappa(\mathcal{D}_2\circ\mathcal{D}_1)\\
            &=\sum_{\beta\in N^{|O_1|}} \sum_{\kappa_1\succ\alpha\cup\beta} \sum_{\kappa_2\succ\beta\cup\gamma} \sgn_{\kappa_1}(\mathcal{D}_1)\sgn_{\kappa_2}(\mathcal{D}_2) \psi_{\kappa_1}(\mathcal{D}_1) \psi_{\kappa_2}(\mathcal{D}_2)\\
            &=\sum_{\beta\in N^{|O_1|}}
                \chi_{\alpha\cup\beta}(\mathcal{D}_1)\chi_{\beta\cup\gamma}(\mathcal{D}_2). \qedhere
            \end{align*}
        \end{proof}
        \end{lemma}
    Returning to the proof of the theorem, since by definition
        $$
        f_{\mathcal{D}_2\circ\mathcal{D}_1}(\bse_{\alpha}) = \sum_{\gamma\in N^{|O_2|}} \chi_{\alpha\cup\gamma}(\mathcal{D}_2\circ\mathcal{D}_1) \bse_{\gamma},
        $$
    it follows from the lemma and \eqref{eq:composition-proof1} that $f_{\mathcal{D}_2}\circ f_{\mathcal{D}_1} = f_{\mathcal{D}_2\circ\mathcal{D}_1}.$
\end{proof}
\end{theorem}
Intuitively, this result means that a trace diagram's function may be understood by breaking the diagram up into little pieces and gluing them back together. For example, the diagram in the introduction is decomposed as follows:
   $$
    \tikz{\node[vertex]at(0,.5){}edge[bend left](.5,0)edge[bend right](-.5,0);
        \node[vertex]at(1.3,.5){}edge[bend right](.8,0)edge[bend left](1.8,0)edge[bend right=90,looseness=1.5](0,.5);}
    =
    \tikz[shift={(0,.1)}]{\draw(0,0)to[out=90,in=90,looseness=2](1,0);}
    \circ\left(
    \tikz{\node[vertex]at(0,.6){}edge(0,1)edge[bend right](-.5,0)edge[bend left](.5,0);}
    \otimes
    \tikz{\node[vertex]at(0,.6){}edge(0,1)edge[bend right](-.5,0)edge[bend left](.5,0);}
    \right).
    $$
This is why the input $\bfu\otimes\bfv\otimes\bfw\otimes\mathbf{x}$ is mapped by the diagram to $(\bfu\times\bfv)\cdot(\bfw\times\mathbf{x})$.

\subsection{Trace Diagram Relations}

\begin{definition}
    A \emph{(framed) trace diagram relation} is a summation $\sum_{\mathcal{D}} c_{\mathcal{D}} \mathcal{D}\in\mathfrak{D}(I,O)$ of framed trace diagrams for which $\sum_{\mathcal{D}} c_{\mathcal{D}} f_{\mathcal{D}} = 0$.
\end{definition}

Under Theorem \ref{thm:categorical}, one can apply trace diagram relations locally on small pieces of larger diagrams. This is exactly what was done in the introduction using the dot and cross product diagrams of Propositions \ref{p:diagram-cap} and \ref{p:diagram-cross}.

Trace diagram relations also exist for \emph{unframed} diagrams, provided the degree-1 vertices are ordered. Let $\mathfrak{D}(m)$ denote the free $\F$-module over tensor diagrams with $m$ ordered degree-1 vertices. Recall that a framing is a partition of these vertices into a set of inputs and a set of outputs. This  provides a mapping $\mathfrak{D}(m)\to\mathfrak{D}(I,O)$ defined whenever $I+O=m$, which we call a \emph{framing}.

\begin{definition}
    A \emph{(general) trace diagram relation} is a summation $\sum_{\mathcal{D}} c_{\mathcal{D}} \mathcal{D}\in\mathfrak{D}(m)$ that restricts under some partition to a framed trace diagram relation.
\end{definition}

\begin{theorem}\label{thm:relations}
    Given a framing $\mathfrak{D}(m)\to\mathfrak{D}(I,O)$, every (general) trace diagram relation in $\mathfrak{D}(m)$ maps to a (framed) trace diagram relation in $\mathfrak{D}(I,O)$.
\begin{proof}
    By Definition \ref{d:diagram-function}, the weights of a function depend only on the leaf labels, and not on the partition or framing of the diagram. Since the weights are the same under different partitions, the relations do not depend on the framing.
\end{proof}
\end{theorem}

The fact that diagrammatic relations are \emph{independent} of framing is very powerful. One may sometimes read off several identities of multilinear algebra from the same diagrammatic relation, as was done in the introduction with \eqref{eq:3binor}. Here is another identity of 3-dimensional vectors:
\begin{example*}
    Using an alternate framing of \eqref{eq:3binor},
    $$
    \tikz[scale=.6]{
        \node[vertex](node)at(.5,1.5){};\node[vertex](node2)at(1,2.5){};
        \draw[](0,0)node[vector]{$\bf u$}to[out=90,in=-135](node);
        \draw[](1,0)node[vector]{$\bf v$}to[out=90,in=-45](node);
        \draw(node)to[out=90,in=-135](node2);
        \draw[](2,0)node[vector]{$\bf w$}to[out=90,in=-45](node2);
        \draw(node2)to[out=90,in=-90](1,3.25);
    }
    =\tikz[scale=.6]{
        \draw[](0,0)node[vector]{$\bf u$}to[out=90,in=90,looseness=2](2,0)node[vector]{$\bf w$};
        \draw[](1,0)node[vector]{$\bf v$}to[out=90,in=-90](1,3.25);
    }
    -\tikz[scale=.6]{
        \draw[](0,0)node[vector]{$\bf u$}to[out=90,in=-90](1,3.25);
        \draw[](1,0)node[vector]{$\bf v$}to[out=90,in=90,looseness=2.5](2,0)node[vector]{$\bf w$};
    }.
    $$
    This proves the following identity:
    $$
    (\bfu\times\bfv)\times\bfw = (\bfu\cdot\bfw) \bfv - (\bfv\cdot\bfw) \bfu.
    $$
\end{example*}

It is even possible for certain diagrams to be decomposed in multiple ways, leading to algebraic identities.
\begin{example*}
    The single diagram
    $$
    \tikz[scale=.6]{
        \node[vertex](node)at(.5,2){};
        \draw[](0,0)node[vector]{$\bf u$}to[out=90,in=-135](node);
        \draw[](1,0)node[vector]{$\bf v$}to[out=90,in=-45](node);
        \draw[](2,0)node[vector]{$\bf w$}to(2,2);\draw(2,2)to[out=90,in=90,looseness=1.5](node);
    }
    =\tikz[scale=.6]{
        \node[vertex](node)at(1.5,2){};
        \draw[](0,0)node[vector]{$\bf u$}to(0,2);\draw(0,2)to[out=90,in=90,looseness=1.5](node);
        \draw[](1,0)node[vector]{$\bf v$}to[out=90,in=-135](node);
        \draw[](2,0)node[vector]{$\bf w$}to[out=90,in=-45](node);
    }
    =\tikz[scale=.6]{
        \node[vertex](node)at(1,2.5){};
        \draw[](0,0)node[vector]{$\bf u$}to[out=90,in=-90](1.5,1.75);\draw(1.5,1.75)to[out=90,in=-45](node);
        \draw[](1,0)node[vector]{$\bf v$}to[out=90,in=-90](2.5,2)to(2.5,2)to[out=90,in=90,looseness=2](node);
        \draw(2,0)node[vector]{$\bf w$}to[out=90,in=-90](.5,1.75);\draw[](.5,1.75)to[out=90,in=-135](node);
    }
    =\tikz[scale=.7]{
        \node[vertex](node)at(1,2){};
        \draw[](0,0)node[vector]{$\bf u$}to[out=90,in=-135](node);
        \draw[](1,0)node[vector]{$\bf v$}to[out=90,in=-90](node);
        \draw[](2,0)node[vector]{$\bf w$}to[out=90,in=-45](node);
    }
    $$
    implies the vector identities
    $$(\bfu\times\bfv)\cdot\bfw=\bfu\cdot(\bfv\times\bfw)=(\bfw\times\bfu)\cdot\bfv=\det[\bfu\:\bfv\:\bfw].$$
    (The fact that
    $\tikz[scale=.6]{
        \node[vertex](node)at(1,1.5){};
        \draw[](0,0)node[vector]{$\bf u$}to[out=90,in=-135](node);
        \draw[](1,0)node[vector]{$\bf v$}to[out=90,in=-90](node);
        \draw[](2,0)node[vector]{$\bf w$}to[out=90,in=-45](node);
    }=\det[\bfu\:\bfv\:\bfw]$
    will be proven in the next section.)
\end{example*}

\ifnum0=1\input{symmetries-unused.tex}\fi

\section{Diagrammatic Building Blocks}\label{s:building-blocks}

This section builds a library of local diagrammatic relations that are needed to reason about general diagrams.

\begin{notation}
    Let $N=\{1,2,\ldots,n\}$.  Given an ordered $k$-tuple $\alpha=(\alpha_1,\alpha_2,\ldots,\alpha_k)\in N^k$ consisting of distinct elements of $N$, let $\overset{\leftarrow}{\alpha}$ denote $(\alpha_k,\ldots,\alpha_2,\alpha_1)$. The switch between $\alpha$ and $\overleftarrow\alpha$ requires $\lfloor\frac{n}{2}\rfloor$ transpositions, where $\lfloor\frac{n}{2}\rfloor=\frac{n}{2}$ if $n$ is even and $\lfloor\frac{n}{2}\rfloor=\frac{n-1}{2}$ if $n$ is odd, and so $\sgn(\overleftarrow\alpha)=\tdsign{n}\sgn(\alpha).$

    Let $S_{\alpha}^c$ represent the set of permutations of $N\setminus\{\alpha_1,\alpha_2,\ldots,\alpha_k\}$. If $\beta=(\beta_1,\beta_2,\ldots,\beta_{n-k})\in S_{\alpha}^c$, let $(\alpha\:\overleftarrow\beta)$ denote the permutation
    $$(\alpha\:\overleftarrow\beta)
    \equiv
    \begin{pmatrix}
        1 & \cdots & k  & k+1 & \cdots & n\\
        \alpha_1 & \cdots & \alpha_k & \beta_{n-k} & \ldots & \beta_1
    \end{pmatrix}.$$
\end{notation}

\begin{proposition}\label{p:one-node-function}
    If $\alpha\in N^k$ has no repeated elements, then
    $$
    \tikz{
        \node[ciliatednode=180](va)at(0,.5){}
            edge[bend left](-.7,1.2)edge[bend right](.4,1.2)edge[bend right](.7,1.2)
            edge[bend right](-.7,-.2)edge[bend right](-.4,-.2)edge[bend left](.7,-.2);
        \draw[dotdotdot](-.6,.9)to node[toplabel]{$n\!-\!k$}(.4,.9);
        \draw[dotdotdot](-.4,.1)to node[bottomlabel]{$k$}(.6,.1);
    }
	:\bse_{\alpha} \longmapsto \sum_{\beta\in S_\alpha^c} \sgn(\alpha\:\overleftarrow\beta)\bse_{\beta}.
    $$
    If $\alpha\in N^k$ has any repeated elements, then the diagram maps $\bse_\alpha$ to 0.
\begin{proof}
    By Definition \ref{d:diagram-function}, the image of $\bse_\alpha$ is automatically 0 if there are repeated elements, since the signature at the node is 0. Otherwise, the diagram maps $\bse_\alpha$ to
    $$
    \sum_{\beta\in N^{n-k}}
            \chi_{\alpha\cup\beta}(\mathcal{D}) \bse_{\beta}
    =\sum_{\beta\in N^{n-k}} \sum_{\kappa\succ\alpha\cup\beta}\sgn_\kappa(\mathcal{D}) \bse_{\beta}
    =\sum_{\beta\in N^{n-k}} \sgn_{\alpha\cup\beta}(\mathcal{D}) \bse_{\beta}
    $$
    Since there are no matrices in the diagram, the coefficient of the coloring is 1. Note that $\alpha\cup\beta$ is a coloring of all edges of the diagram. If $\beta$ includes any of the same elements as $\alpha$, the signature of the coloring is zero. Therefore, we may restrict to the summation in which $\beta\in S_\alpha^c$. In this situation, $\alpha\cup\beta$ is a proper coloring of the entire diagram, and the signature is then
    $$
    \sgn_{\alpha\cup\beta}(\mathcal{D})=\sgn(\alpha\:\overset\leftarrow\beta). \qedhere
    $$
\end{proof}
\end{proposition}

Some special cases of this result are particularly useful. When $k=n$, this proposition states that
    \begin{equation}\label{eq:n-cap}
    \tikz[shift={(0,.3)}]{
        \node[ciliatednode=160](va)at(0,.5){}
            edge[bend right](-.7,-.2)edge[bend right](-.4,-.2)edge[bend left](.7,-.2);
        \draw[dotdotdot](-.4,.1)to node[bottomlabel]{$n$}(.6,.1);
    }
    :\bse_\alpha \longmapsto \sgn(\alpha) = \det(\bse_{\alpha_1} \cdots \bse_{\alpha_n}).
    \end{equation}
Therefore, by linear extension:
    \begin{equation}\label{eq:n-cap2}
    \tikz[scale=1.5,shift={(0,.2)}]{
        \node[ciliatednode=160](va)at(0,.5){}
            edge[bend right]node[small vector,pos=1]{$\bfu_1$}(-.7,-.2)
            edge[bend right]node[small vector,pos=1]{$\bfu_2$}(-.4,-.2)
            edge[bend left]node[small vector,pos=1]{$\bfu_n$}(.7,-.2);
        \draw[dotdotdot](-.4,.1)to node[bottomlabel]{$n$}(.6,.1);
    }
    = \det(\bfu_{1} \cdots \bfu_{n}).
    \end{equation}
When $k=0$, Proposition \ref{p:one-node-function} states that
    \begin{equation}\label{eq:n-cup}
    \tikz[shift={(0,-.3)}]{
        \node[ciliatednode=200](va)at(0,.5){}
            edge[bend left](-.7,1.2)edge[bend left](-.4,1.2)edge[bend right](.7,1.2);
        \draw[dotdotdot](-.4,.9)to node[toplabel]{$n$}(.6,.9);
    }
    :1 \longmapsto \sum_{\beta\in S_n} \sgn(\overset\leftarrow\beta) \bse_{\beta}
    = (-1)^{\lfloor\frac{n}{2}\rfloor} \sum_{\beta\in S_n} \sgn(\beta) \bse_{\beta}.
    \end{equation}
Also, the case $k=n-1$ provides a generalization of the three-dimensional cross product.

\begin{proposition}\label{p:two-node-function}
    If $\alpha\in N^k$ has no repeated elements, then
    $$
    \tikz[heighttwo]{
        \node[ciliatednode=180](vb)at(0,1.5){}
            edge[bend left](-.7,2.2)edge[bend left](-.4,2.2)edge[bend right](.7,2.2);
        \node[ciliatednode=180](va)at(0,.5){}
            edge[out=160,in=200,looseness=2](vb)edge[out=20,in=-20,looseness=2](vb)edge[out=20,in=-20,looseness=1.2](vb)
            edge[bend right](-.7,-.2)edge[bend right](-.4,-.2)edge[bend left](.7,-.2);
        \draw[dotdotdot](-.4,1.9)to node[toplabel]{$k$}(.6,1.9);
        \draw[dotdotdot](-.6,.8)to node[toplabel,scale=.9,xscale=.9]{$n\!-\!k$}(.4,.8);
        \draw[dotdotdot](-.4,.1)to node[bottomlabel]{$k$}(.6,.1);
    }:\bse_\alpha \longmapsto (-1)^{\lfloor\frac{n}{2}\rfloor}(n-k)!
    \sum_{\sigma\in S_\alpha} \sgn\qptmx{\alpha}{\sigma}\bse_{\sigma(\alpha)},
    $$
    where $\sgn\qptmx{\alpha}{\sigma}=(-1)^t$ when $t$ transpositions are required to transform $\alpha$ into $\sigma$.
    If $\alpha\in N^k$ has any repeated elements, then the diagram maps $\bse_\alpha$ to 0.
\begin{proof}
    Applying Proposition \ref{p:one-node-function} twice (and noting that if $\beta\in S_\alpha^c then S_\beta^c=S_\alpha$), the image of $\bse_\alpha$ is
        $$\sum_{\beta\in S_\alpha^c} \sgn(\alpha\:\overset\leftarrow\beta) \sum_{\sigma\in S_\alpha} \sgn(\beta\:\overset\leftarrow\sigma) \bse_\sigma.$$

    We claim that
        $\sgn(\alpha\:\overset\leftarrow\beta)\sgn(\beta\:\overset\leftarrow\sigma)
            =\sgn(\alpha\:\overset\leftarrow{\beta'})\sgn(\beta'\:\overset\leftarrow\sigma)$
     for any $\beta,\beta'\in S_\alpha^c$. To see this, consider the process of transposing elements to change $\beta$ into $\beta'$. If this process requires $t$ transpositions, then $\sgn(\beta)=(-1)^t\sgn(\beta')$, which implies both $\sgn(\alpha\:\overset\leftarrow\beta)=(-1)^t\sgn(\alpha\:\overset\leftarrow{\beta'})$ and $\sgn(\beta\:\overset\leftarrow\sigma)=(-1)^t\sgn(\beta'\:\overset\leftarrow\sigma)$. The claim follows.

    Given this claim, every $\beta\in S_\alpha^c$ makes the same contribution to the sum, and the expression reduces to
        $$(n-k)!\sum_{\sigma\in S_\alpha}
            \sgn(\alpha\:\overset\leftarrow\beta)\sgn(\beta\:\overset\leftarrow\sigma)
            \bse_\sigma,$$
    where $\beta$ is an arbitrary element of $S_\alpha^c$. The signature term simplifies as follows:
    \begin{multline*}
    \sgn(\alpha\:\overset\leftarrow\beta)\sgn(\beta\:\overset\leftarrow\sigma)
        =\sgn(\alpha\:\overset\leftarrow\beta)(-1)^{\lfloor\frac{n}{2}\rfloor}\sgn(\sigma\:\overset\leftarrow\beta)\\
        =(-1)^{\lfloor\frac{n}{2}\rfloor}\sgn\qptmx{\alpha}{\sigma}\sgn(\alpha\:\overset\leftarrow\beta)^2
        =(-1)^{\lfloor\frac{n}{2}\rfloor}\sgn\qptmx{\alpha}{\sigma}.\qedhere
    \end{multline*}
\end{proof}
\end{proposition}

The following result depends on the previous proof, and is used repeatedly in later sections.
\begin{lemma}[cut-and-paste lemma]\label{l:cut-paste}
    If $\alpha\in N^k$ has no repeated elements, $\beta\in S_\alpha^c$, and $A$ is any $n\times n$ matrix, then
    \begin{equation}\label{eq:cut-paste}
    \tikz[heighttwo,scale=1.2]{
        \node[ciliatednode=180](vb)at(0,1.6){}
            edge[bend left](-.7,2.2)
            edge[bend left](-.4,2.2)
            edge[bend right](.7,2.2);
        \node[ciliatednode=180](va)at(0,.3){}
            edge[out=160,in=200,looseness=2]node[matrix on edge]{$A$}(vb)
            edge[out=20,in=-20,looseness=2]node[matrix on edge]{$A$}(vb)
            edge[out=20,in=-20,looseness=1.2]node[matrix on edge]{$A$}(vb)
            edge[bend right](-.7,-.2)
            edge[bend right](-.4,-.2)
            edge[bend left](.7,-.2);
        \node[bottomlabel]at(0,-.2){$\alpha_1\alpha_2\cdots\alpha_k$};
        \draw[dotdotdot](-.4,1.9)to node[toplabel]{$k$}(.6,1.9);
        \draw[dotdotdot](-.6,1)to node[toplabel,scale=.8]{$n\!-\!k$}(.4,1);
    }
    = \sgn(\alpha \:\overset\leftarrow\beta)(n-k)!
    \tikz[shift={(0,.6)}]{
        \node[ciliatednode=180](va)at(0,.5){}
            edge[bend left](-.7,1.2)
            edge[bend left](-.4,1.2)
            edge[bend right](.7,1.2)
            edge[bend right]node[matrix on edge,pos=.7]{$A$}(-.7,-.5)
            edge[bend left]node[reverse matrix on edge,pos=.7]{$A$}(.4,-.5)
            edge[bend left]node[reverse matrix on edge,pos=.7]{$A$}(.7,-.5);
        \node[bottomlabel]at(.1,-.5){$\beta_1\:\cdots\:\beta_{n-k}$};
        \draw[dotdotdot](-.4,.9)to node[toplabel]{$k$}(.6,.9);
    }.
    \end{equation}
    If $\alpha\in N^k$ has repeated elements, then the diagram maps $\bse_\alpha$ to 0.
\begin{proof}
    By Proposition \ref{p:one-node-function}, the lefthand side of \eqref{eq:cut-paste} evaluates to
    $$\sum_{\beta\in S_\alpha^c} \sgn(\alpha \:\overset\leftarrow\beta)
    \tikz{
        \node[ciliatednode=180](va)at(0,.5){}
            edge[bend left](-.7,1.2)
            edge[bend left](-.4,1.2)
            edge[bend right](.7,1.2)
            edge[bend right]node[matrix on edge,pos=.7]{$A$}(-.7,-.5)
            edge[bend left]node[reverse matrix on edge,pos=.7]{$A$}(.4,-.5)
            edge[bend left]node[reverse matrix on edge,pos=.7]{$A$}(.7,-.5);
        \node[bottomlabel]at(.1,-.5){$\beta_1\:\cdots\:\beta_{n-k}$};
        \draw[dotdotdot](-.4,.9)to node[toplabel]{$k$}(.6,.9);
    }.$$
    As in the proof of Proposition \ref{p:two-node-function}, the result is true because every choice of $\beta$ contributes the same value to the summand. In this case,
    a transposition of elements of $\beta$ corresponds to swapping two of the strands labeled by $\beta_i$ in the diagram. But swapping two strands in the diagram leads to a change of signature at the node. In particular, if $\beta,\beta'\in S_\alpha^c$ are related by $t$ transpositions, then $\sgn(\alpha\:\overset\leftarrow\beta)=(-1)^t\sgn(\alpha\:\overset\leftarrow{\beta'})$ and
    $$
    \tikz{
        \node[ciliatednode=180](va)at(0,.5){}
            edge[bend left](-.7,1.2)
            edge[bend left](-.4,1.2)
            edge[bend right](.7,1.2)
            edge[bend right]node[matrix on edge,pos=.7]{$A$}(-.7,-.5)
            edge[bend left]node[reverse matrix on edge,pos=.7]{$A$}(.4,-.5)
            edge[bend left]node[reverse matrix on edge,pos=.7]{$A$}(.7,-.5);
        \node[bottomlabel]at(.1,-.5){$\beta_1\:\cdots\:\beta_{n-k}$};
        \draw[dotdotdot](-.4,.9)to node[toplabel]{$k$}(.6,.9);
    }=(-1)^t
     \tikz{
        \node[ciliatednode=180](va)at(0,.5){}
            edge[bend left](-.7,1.2)
            edge[bend left](-.4,1.2)
            edge[bend right](.7,1.2)
            edge[bend right]node[matrix on edge,pos=.7]{$A$}(-.7,-.5)
            edge[bend left]node[reverse matrix on edge,pos=.7]{$A$}(.4,-.5)
            edge[bend left]node[reverse matrix on edge,pos=.7]{$A$}(.7,-.5);
        \node[bottomlabel]at(.1,-.5){$\beta'_1\:\cdots\:\beta'_{n-k}$};
        \draw[dotdotdot](-.4,.9)to node[toplabel]{$k$}(.6,.9);
    }.$$
    Consequently, the summation may be replaced by the number of elements in $S_\alpha^c$, which is $(n-k)!$.
\end{proof}
\end{lemma}
This result is called the ``cut-and-paste lemma'' because it allows nodes to be removed or added on to certain parts of a trace diagram. It will be used frequently in later sections.

The following result is vital to manipulating matrices within diagrams. Note that both statements in the theorem are \emph{general} trace diagram relations.
\begin{proposition}[matrix action at nodes]\label{p:node-equivariance}
    If $A$ is any $n\times n$ matrix, then
    \begin{equation}\label{eq:node-equivariance}
    \tikz[scale=1.5,shift={(0,.25)}]{
        \node[ciliatednode=180](va)at(0,.5){}
            edge[bend right]node[matrix on edge,pos=.7]{$A$}(-.7,-.2)
            edge[bend right]node[matrix on edge,pos=.7]{$A$}(-.4,-.2)
            edge[bend left]node[reverse matrix on edge,pos=.7]{$A$}(.7,-.2);
        \draw[dotdotdot](-.4,.1)to node[bottomlabel]{$n$}(.6,.1);
        \draw[dotdotdot](-.4,.1)to(.6,.1);
    }=\det(A)
    \tikz[scale=1.5,shift={(0,.25)}]{
        \node[ciliatednode=180](va)at(0,.5){}
            edge[bend right](-.7,-.2)edge[bend right](-.4,-.2)edge[bend left](.7,-.2);
        \draw[dotdotdot](-.4,.1)to node[bottomlabel]{$n$}(.6,.1);
        \draw[dotdotdot](-.4,.1)to(.6,.1);
    }.
    \end{equation}
    If $A$ is an invertible $n\times n$ matrix, and
        \tikz{\draw(0,0)to node[matrix on edge]{$\overline{A}$}(0,1);}
    represents its inverse $A^{-1}$, then
    \begin{equation}\label{eq:node-equivariance-inverse}
    \tikz[scale=1.3]{
        \node[ciliatednode=180](va)at(0,.5){}
            edge[bend left](-.7,1.2)
            edge[bend left](-.4,1.2)
            edge[bend right](.7,1.2)
            edge[bend right]node[matrix on edge,pos=.7]{$A$}(-.7,-.5)
            edge[bend left]node[reverse matrix on edge,pos=.7]{$A$}(.4,-.5)
            edge[bend left]node[reverse matrix on edge,pos=.7]{$A$}(.7,-.5);
        \draw[dotdotdot](-.4,.9)to node[toplabel]{$k$}(.6,.9);
        \draw[dotdotdot](-.6,-.4)to node[bottomlabel,scale=.8]{$n\!-\!k$}(.4,-.4);
    }=\det(A)
    \tikz[scale=1.3]{
        \node[ciliatednode=180](va)at(0,.2){}
            edge[bend left]node[matrix on edge,pos=.65]{$\overline{A}$}(-.7,1.2)
            edge[bend left]node[matrix on edge,pos=.65]{$\overline{A}$}(-.4,1.2)
            edge[bend right]node[reverse matrix on edge,pos=.65]{$\overline{A}$}(.7,1.2)
            edge[bend right](-.7,-.5)
            edge[bend left](.4,-.5)
            edge[bend left](.7,-.5);
        \draw[dotdotdot](-.6,-.2)to node[bottomlabel,scale=.8]{$n\!-\!k$}(.4,-.2);
        \draw[dotdotdot](-.4,.9)to node[toplabel]{$k$}(.6,.9);
    }.
    \end{equation}
\begin{proof}
    Theorem \ref{thm:categorical} greatly simplifies this proof, since it allows one to compute a diagram's function by starting from an arbitrary input at the bottom, and working upward through the diagram. Let $\bse_\alpha\in N^n$ represent a basis input to \eqref{eq:node-equivariance} and let $A_i$ denotes the $i$th column of $A$. Then
    \begin{equation*}\label{eq:node-step1}
    \tikz[scale=1.8,shift={(0,.2)}]{
        \node[ciliatednode=180](va)at(0,.5){}
            edge[bend right]node[matrix on edge]{$A$}(-.7,-.2)
            edge[bend right]node[matrix on edge]{$A$}(-.4,-.2)
            edge[bend left]node[reverse matrix on edge]{$A$}(.7,-.2);
        \node[bottomlabel]at(0,-.2){$\alpha_1\:\alpha_2\:\:\:\cdots\:\:\:\alpha_n$};
    }=\tikz[scale=1.8,shift={(0,.2)}]{
        \node[ciliatednode=180](va)at(0,.5){}
            edge[bend right]node[small vector,pos=1]{$A_{\alpha_1}$}(-.7,-.2)
            edge[bend right]node[small vector,pos=1]{$A_{\alpha_2}$}(-.4,-.2)
            edge[bend left]node[small vector,pos=1]{$A_{\alpha_n}$}(.7,-.2);
        \draw[dotdotdot](-.3,-.2)to(.6,-.2);
    }=\det(A_{\alpha_1}A_{\alpha_2}\cdots A_{\alpha_n}),
    \end{equation*}
    where the last step follows from \eqref{eq:n-cap2}. Observe that
    $$\det(A_{\alpha_1}A_{\alpha_2}\cdots A_{\alpha_n})=\sgn(\alpha)\det(A),$$
    since the number of transpositions required to restore $\alpha$ to the identity permutation is the same number of column switches required to restore the matrix $(A_{\alpha_1}A_{\alpha_2}\cdots A_{\alpha_n})$ to the original matrix $A$. The proof is completed by noting that $\sgn(\alpha)\det(A)$ is the value of the righthand side of \eqref{eq:node-equivariance} for the input $\bse_\alpha$.

    The second statement \eqref{eq:node-equivariance-inverse} follows from the first, by inserting an explicit copy of the identity matrix in the form of $A A^{-1}$ on the top strands, and then applying \eqref{eq:node-equivariance}:
    $$\tikz[scale=1.3]{
        \node[ciliatednode=180](va)at(0,.5){}
            edge[bend left](-.7,1.2)
            edge[bend left](-.4,1.2)
            edge[bend right](.7,1.2)
            edge[bend right]node[matrix on edge,pos=.7]{$A$}(-.7,-.5)
            edge[bend left]node[reverse matrix on edge,pos=.7]{$A$}(.4,-.5)
            edge[bend left]node[reverse matrix on edge,pos=.7]{$A$}(.7,-.5);
        \draw[dotdotdot](-.4,.9)to node[toplabel]{$k$}(.6,.9);
        \draw[dotdotdot](-.6,-.4)to node[bottomlabel,scale=.8]{$n\!-\!k$}(.4,-.4);
    }=
    \tikz[scale=1.7]{
        \node[ciliatednode=180](va)at(0,.2){}
            edge[bend left]node[matrix on edge,pos=.3,scale=.9]{$A$}node[matrix on edge,pos=.7,scale=.9]{$\overline{A}$}(-.7,1.2)
            edge[bend left]node[matrix on edge,pos=.3,scale=.9]{$A$}node[matrix on edge,pos=.7,scale=.9]{$\overline{A}$}(-.4,1.2)
            edge[bend right]node[reverse matrix on edge,pos=.3,scale=.9]{$A$}node[reverse matrix on edge,pos=.7,scale=.9]{$\overline{A}$}(.7,1.2)
            edge[bend right]node[matrix on edge,pos=.6,scale=.9]{$A$}(-.7,-.5)
            edge[bend left]node[reverse matrix on edge,pos=.6,scale=.9]{$A$}(.4,-.5)
            edge[bend left]node[reverse matrix on edge,pos=.6,scale=.9]{$A$}(.7,-.5);
        \draw[dotdotdot](-.4,.9)to node[toplabel]{$k$}(.6,.9);
        \draw[dotdotdot](-.6,-.4)to node[bottomlabel,scale=.8]{$n\!-\!k$}(.4,-.4);
    }=\det(A)
    \tikz[scale=1.3]{
        \node[ciliatednode=180](va)at(0,.2){}
            edge[bend left]node[matrix on edge,pos=.65]{$\overline{A}$}(-.7,1.2)
            edge[bend left]node[matrix on edge,pos=.65]{$\overline{A}$}(-.4,1.2)
            edge[bend right]node[reverse matrix on edge,pos=.65]{$\overline{A}$}(.7,1.2)
            edge[bend right](-.7,-.5)
            edge[bend left](.4,-.5)
            edge[bend left](.7,-.5);
        \draw[dotdotdot](-.6,-.2)to node[bottomlabel,scale=.8]{$n\!-\!k$}(.4,-.2);
        \draw[dotdotdot](-.4,.9)to node[toplabel]{$k$}(.6,.9);
    }. \qedhere
    $$
\end{proof}
\end{proposition}

\begin{example*}
    One can use \eqref{eq:node-equivariance} to prove that $\det(AB)=\det(A)\det(B)$. Applying \eqref{eq:node-equivariance} directly gives
    $$
    \tikz[scale=1.8,shift={(0,.1)}]{
        \node[ciliatednode=180](va)at(0,.6){}
            edge[bend right]node[matrix on edge,pos=.7,scale=.9]{$AB$}(-.7,-.2)
            edge[bend right]node[matrix on edge,pos=.7,scale=.9]{$AB$}(-.4,-.2)
            edge[bend left]node[reverse matrix on edge,pos=.7,scale=.9]{$AB$}(.7,-.2);
        \draw[dotdotdot](-.4,.1)to node[bottomlabel]{$n$}(.6,.1);
        \draw[dotdotdot](-.4,.1)to(.6,.1);
    }=\det(AB)
    \tikz[scale=1.8,shift={(0,.1)}]{
        \node[ciliatednode=180](va)at(0,.5){}
            edge[bend right](-.7,-.2)edge[bend right](-.4,-.2)edge[bend left](.7,-.2);
        \draw[dotdotdot](-.4,.1)to node[bottomlabel]{$n$}(.6,.1);
        \draw[dotdotdot](-.4,.1)to(.6,.1);
    }.$$
    On the other hand, one may use the fact that
    $\tikz[heightoneonehalf]{
            \draw(0,0)to(0,1.5)
                node[matrix on edge,pos=.25]{$B$}
                node[matrix on edge,pos=.75]{$A$};}
        =
        \tikz{\draw(0,0)to[wavyup](.05,1)node[matrix on edge,pos=.5]{$AB$};}$
    to write the same diagram as
    $$\tikz[scale=1.5,shift={(0,.25)}]{
        \node[ciliatednode=180](va)at(0,.6){}
            edge[bend right]node[matrix on edge,pos=.7,scale=.7]{$AB$}(-.7,-.2)
            edge[bend right]node[matrix on edge,pos=.7,scale=.7]{$AB$}(-.4,-.2)
            edge[bend left]node[reverse matrix on edge,pos=.7,scale=.7]{$AB$}(.7,-.2);
        \draw[dotdotdot](-.4,.1)to node[bottomlabel]{$n$}(.6,.1);
        \draw[dotdotdot](-.4,.1)to(.6,.1);
    }=\tikz[scale=1.5,shift={(0,.25)}]{
        \node[ciliatednode=180](va)at(0,.7){}
            edge[bend right]node[matrix on edge,pos=.4,scale=.8]{$A$}node[matrix on edge,pos=.8,scale=.8]{$B$}(-.7,-.2)
            edge[bend right]node[matrix on edge,pos=.4,scale=.8]{$A$}node[matrix on edge,pos=.8,scale=.8]{$B$}(-.4,-.2)
            edge[bend left]node[reverse matrix on edge,pos=.4,scale=.8]{$A$}node[reverse matrix on edge,pos=.8,scale=.8]{$B$}(.7,-.2);
        \draw[dotdotdot](-.4,.1)to node[bottomlabel]{$n$}(.6,.1);
        \draw[dotdotdot](-.4,.1)to(.6,.1);
    }
    =\det(A)
    \tikz[scale=1.5,shift={(0,.25)}]{
        \node[ciliatednode=180](va)at(0,.6){}
            edge[bend right]node[matrix on edge,pos=.7,scale=.8]{$B$}(-.7,-.2)
            edge[bend right]node[matrix on edge,pos=.7,scale=.8]{$B$}(-.4,-.2)
            edge[bend left]node[reverse matrix on edge,pos=.7,scale=.8]{$B$}(.7,-.2);
        \draw[dotdotdot](-.4,.1)to node[bottomlabel]{$n$}(.6,.1);
        \draw[dotdotdot](-.4,.1)to(.6,.1);
    }
    =\det(A)\det(B)
    \tikz[scale=1.5,shift={(0,.25)}]{
        \node[ciliatednode=180](va)at(0,.5){}
            edge[bend right](-.7,-.2)edge[bend right](-.4,-.2)edge[bend left](.7,-.2);
        \draw[dotdotdot](-.4,.1)to node[bottomlabel]{$n$}(.6,.1);
        \draw[dotdotdot](-.4,.1)to(.6,.1);
    }.$$
\end{example*}
One can similarly apply Proposition \ref{p:transpose-diagram} to the relation \eqref{eq:node-equivariance} to show that $\det(A^T)=\det(A)$.

\begin{proposition}[Determinant Diagram]\label{p:determinant-diagram}
\begin{equation}\label{detdiagram}
    \tikz[heightoneonehalf]{
            \node[ciliatednode=170](topnode)at(0,1.5){};
            \node[ciliatednode=-170](bottomnode)at(0,0){};
            \draw(bottomnode)to[bend left](-.7,.75)node[small matrix]{$A$}to[bend left](topnode);
            \draw(bottomnode)to[bend left](-.4,.75)node[small matrix]{$A$}to[bend left](topnode);
            \draw(bottomnode)to[bend right](.7,.75)node[small matrix]{$A$}to[bend right](topnode);
            \draw[dotdotdot](-.3,.75)to node[toplabel]{$n$}(.6,.75);
        }
    =(-1)^{\lfloor\frac{n}{2}\rfloor} n! \det(A).
\end{equation}
\begin{proof}
    Proposition \ref{p:node-equivariance} gives the factor $\det(A)$, while Proposition \ref{p:two-node-function} with $k=0$ gives the factor $(-1)^{\lfloor\frac{n}{2}\rfloor} n!$.
\end{proof}
\end{proposition}

\section{Matrix Minors}\label{s:minors}

This section reveals the fundamental role of matrix minors in trace diagram functions. We begin with notation and a review of matrix minors. A good classical treatment of matrix minors is given in Section 2.4 of \cite{LT85}.

\subsection{Matrix Minors and Cofactors}
Let $A$ be an $n\times n$ matrix over a field $\F$ with
	$$A = \begin{pmatrix}
		a_{11} & a_{12} & \cdots & a_{1n}\\
		a_{21} & a_{22} & \cdots & a_{2n}\\
		\vdots & \vdots & \ddots & \vdots\\
		a_{n1} & a_{n2} & \cdots & a_{nn}
		\end{pmatrix}.$$
A \emph{submatrix} of a matrix $A$ is a smaller matrix formed by ``crossing out'' a number of rows and columns in $A$.

Let $N \equiv \{1,2,\ldots,n\}$. Let $I=(I_1,\ldots,I_{k_1})$ and $J=(J_1,\ldots,J_{k_2})$ be ordered subsets of $N$ in which $1\leq I_1<\cdots<I_{k_1}\leq n$ and similarly for $J$. Let $A_{I,J}$ denote the submatrix formed from the rows in $I$ and the columns in $J$. The \emph{complementary submatrix} $A^c_{I,J}$ is formed by crossing out the rows in $I$ and the columns in $J$. For $n\ge3$, the \emph{interior} $\mathrm{int}(A)$ is the submatrix $A^c_{(1,n),(1,n)}$.

\begin{example*}
    Let $I=(1,2)$ and $J=(3,4)$. If
    \begin{align*}
	   A=\begin{pmatrix}
		  a&b&c&d\\
    		e&f&g&h\\
	   	i&j&k&l\\
		  m&n&o&p
	   \end{pmatrix},
    \end{align*}
    then $A_{I,J}=\qptmx{c&d}{g&h}$, $A^c_{I,J}=\qptmx{i&j}{m&n}$, and $\mathrm{int}(A)=\qptmx{f&g}{j&k}$.
\end{example*}

\begin{definition}\label{d:matrix-minor}
    If $I$ and $J$ have the same number of entries, the \emph{minor} $[A_{I,J}]$ is the determinant of the submatrix $A_{I,J}$. The \emph{complementary minor} $[A^c_{I,J}]$ is the determinant of the complementary submatrix $A^c_{I,J}$.
\end{definition}
A direct formula for the $k\times k$ minor is
    \begin{equation}\label{eq:matrix-minor}
    [A_{I,J}] = \sum_{\sigma\in S_k} \mathrm{sgn}(\sigma) a_{I_1,J_{\sigma(1)}} a_{I_2,J_{\sigma(2)}} \cdots a_{I_k,J_{\sigma(k)}}.
    \end{equation}
In the above example, $[A_{I,J}]=ch-gd$.

\begin{definition}\label{d:cofactor}
    The \emph{$(i,j)$-cofactor} of $A$ is
	   $$C_{ij} \equiv (-1)^{i+j} [A^c_{i,j}].$$
    The \emph{$(I,J)$-cofactor} of $A$ is
        $$C_{I,J}=(-1)^{I_1+\cdots+I_k+J_1+\cdots+J_k}[A^c_{I,J}].$$
    The \emph{adjugate} (or \emph{adjoint}) $\adj(A)$ of a square matrix is the matrix comprised of entries $(\adj(A))_{ij}\equiv C_{ji}$.
\end{definition}
A student often sees cofactors first in the \emph{cofactor expansion} formula useful for by-hand calculations of the determinant:
	\begin{equation}\label{eq:laplace-expansion}
	\det(A)=\sum_{j=1}^n a_{ij} C_{ij},
	\end{equation}
where $i\in N$ is an arbitrary row. Adjugates are sometimes used to compute the matrix inverse since $A^{-1} = \frac{1}{\det(A)}\adj(A)$ when $A$ is invertible.

\subsection{Diagrams for Matrix Minors}

\begin{proposition}\label{p:minor-diagram}
    Let $A$ be an $n\times n$ matrix. Then
    \begin{equation}\label{eq:minor-diagram}
    [A_{I,J}]
    =\sgn(J^c\:\overset\leftarrow{J})
    \tikz{
        \node[ciliatednode=180](va)at(0,.2){}
            edge[bend left]node[reverse matrix on edge,pos=.65]{$A$}(-.7,1.2)
            edge[bend left]node[reverse matrix on edge,pos=.65]{$A$}(-.4,1.2)
            edge[bend right]node[matrix on edge,pos=.65]{$A$}(.7,1.2)
            \foreach\xa/\xb in{right/-.7,left/.4,left/.7}{edge[bend \xa](\xb,-.5)};
        \draw(0,1.2)node[toplabel]{$I_1I_2\cdots I_k$};
        \draw(.2,-.5)node[bottomlabel]{$J_1^c\:\cdots\:J_{n-k}^c$};}
    =\sgn(I^c\:\overset\leftarrow{I})
    \tikz{
        \node[ciliatednode=180](va)at(0,.2){}
            edge[bend left]node[matrix on edge,pos=.65]{$A$}(-.7,1.2)
            edge[bend left]node[matrix on edge,pos=.65]{$A$}(-.4,1.2)
            edge[bend right]node[reverse matrix on edge,pos=.65]{$A$}(.7,1.2)
            edge[bend right](-.7,-.5)
            edge[bend left](.4,-.5)
            edge[bend left](.7,-.5);
        \draw(0,1.2)node[toplabel]{$J_1J_2\cdots J_k$};
        \draw(.2,-.5)node[bottomlabel]{$I_1^c\:\cdots\:I_{n-k}^c$};
    }.
    \end{equation}
\begin{proof}
    By Proposition \ref{p:two-node-function} and the minor formula \eqref{eq:matrix-minor},
    \begin{equation*}
    \tikz[heighttwo]{
        \node[ciliatednode=180](va)at(0,1.3){}
            edge[bend left]node[reverse matrix on edge,pos=.65,scale=.8]{$A$}(-.7,2)
            edge[bend left]node[reverse matrix on edge,pos=.65,scale=.8]{$A$}(-.4,2)
            edge[bend right]node[matrix on edge,pos=.65,scale=.8]{$A$}(.7,2);
        \node[ciliatednode=180](vb)at(0,.2){}
            edge[bend right](-.7,-.5)
            edge[bend right](-.4,-.5)
            edge[bend left](.7,-.5)
            edge[out=160,in=200,looseness=2](va)
            edge[out=20,in=-20,looseness=2](va)
            edge[out=20,in=-20,looseness=1.2](va);
        \draw(0,2)node[toplabel]{$I_1I_2\cdots I_k$};
        \draw[dotdotdot](-.45,.5)to(.3,.5)node[toplabel,scale=.8,pos=.5]{$n\!-\!k$};
        \draw(0,-.5)node[bottomlabel]{$J_1J_2\cdots J_k$};
    }
    = \tdsign{n}(n-k)!\sum_{\sigma\in S_k}\sgn(\sigma)
        \tikz[heighttwo]{
            \foreach\xa/\xb in{-1,-.5,1}{\draw(\xa,0)to(\xa,2)node[matrix on edge,pos=.5]{$A$};}
            \draw[dotdotdot](-.5,1)to(1,1);
            \draw(0,0)node[bottomlabel]{$\:\:J_{\sigma(1)}\:\:\cdots J_{\sigma(k)}$};
            \draw(0,2)node[toplabel]{$I_1\:I_2\:\:\cdots\:\:I_k$};
        }
    = \tdsign{n}(n-k)![A_{I,J}].
    \end{equation*}
    Using the cut-and-paste lemma \eqref{eq:cut-paste}, the same diagram reduces to
    \begin{equation*}
    (n-k)! \sgn(J\:\overset\leftarrow{J^c})
    \tikz{
        \node[ciliatednode=180](va)at(0,.2){}
            edge[bend left]node[reverse matrix on edge,pos=.65]{$A$}(-.7,1.2)
            edge[bend left]node[reverse matrix on edge,pos=.65]{$A$}(-.4,1.2)
            edge[bend right]node[matrix on edge,pos=.65]{$A$}(.7,1.2)
            edge[bend right](-.7,-.5)
            edge[bend left](.4,-.5)
            edge[bend left](.7,-.5);
        \draw(0,1.2)node[toplabel]{$I_1I_2\cdots I_k$};
        \draw(.2,-.5)node[bottomlabel]{$J_1^c\:\cdots\:J_{n-k}^c$};
    }
    =(n-k)! \tdsign{n} \sgn(J^c\:\overset\leftarrow{J})
    \tikz{
        \node[ciliatednode=180](va)at(0,.2){}
            edge[bend left]node[reverse matrix on edge,pos=.65]{$A$}(-.7,1.2)
            edge[bend left]node[reverse matrix on edge,pos=.65]{$A$}(-.4,1.2)
            edge[bend right]node[matrix on edge,pos=.65]{$A$}(.7,1.2)
            edge[bend right](-.7,-.5)
            edge[bend left](.4,-.5)
            edge[bend left](.7,-.5);
        \draw(0,1.2)node[toplabel]{$I_1I_2\cdots I_k$};
        \draw(.2,-.5)node[bottomlabel]{$J_1^c\:\cdots\:J_{n-k}^c$};
    }.
    \end{equation*}
    This verifies the first function. The second case is similar.
\end{proof}
\end{proposition}

The next section requires understanding the following diagrams for the cofactor and the adjugate:
\begin{proposition}\label{p:cofactor-diagram}
    Let $A$ be an $n\times n$ matrix. Then
    \begin{equation}\label{eq:cofactor-diagram}
    C_{I,J}=
    \frac{\tdsign{n}}{(n-k)!}
    \tikz[heighttwo,scale=1.2]{
        \node[ciliatednode=180](vb)at(0,1.6){}
            edge[bend left](-.7,2.2)
            edge[bend left](-.4,2.2)
            edge[bend right](.7,2.2);
        \node[ciliatednode=180](va)at(0,.3){}
            edge[out=160,in=200,looseness=2]node[matrix on edge]{$A$}(vb)
            edge[out=20,in=-20,looseness=2]node[matrix on edge]{$A$}(vb)
            edge[out=20,in=-20,looseness=1.2]node[matrix on edge]{$A$}(vb)
            edge[bend right](-.7,-.2)
            edge[bend right](-.4,-.2)
            edge[bend left](.7,-.2);
        \node[bottomlabel]at(0,-.2){$J_1J_2\cdots J_k$};
        \node[toplabel]at(0,2.2){$I_1I_2\cdots I_k$};
        \draw[dotdotdot](-.6,1)to node[toplabel,scale=.8]{$n\!-\!k$}(.4,1);
    }
    \qquad\text{and}\qquad
    \mathsf{adj}(A)
    =\frac{\tdsign{n}}{(n-1)!}
    \tikz[heighttwo,scale=1.2]{
        \node[ciliatednode=160](vb)at(0,1.6){}edge(0,2.2);
        \node[ciliatednode=200](va)at(0,.3){}
            edge[out=160,in=200,looseness=1.7]node[reverse matrix on edge]{$A$}(vb)
            edge[out=20,in=-20,looseness=1.6]node[reverse matrix on edge]{$A$}(vb)
            edge[out=20,in=-20,looseness=1.1]node[reverse matrix on edge]{$A$}(vb)
            edge(0,-.2);
        \draw[dotdotdot](-.4,.9)to node[toplabel,scale=.7]{$n\!-\!1$}(.2,.9);
    }.
    \end{equation}
\begin{proof}
    By Proposition \ref{p:minor-diagram} and the cut-and-paste lemma \eqref{eq:cut-paste} (and replacing $I$ with $I^c$ and $J$ with $J^c$), the complementary minor is:
    \begin{equation}\label{eq:minor2}
    [A^c_{I,J}]=\sgn(J\:\overset\leftarrow{J^c})
    \tikz{
        \node[ciliatednode=180](va)at(0,.2){}
            edge[bend left]node[reverse matrix on edge,pos=.65]{$A$}(-.7,1.2)
            edge[bend right]node[matrix on edge,pos=.65]{$A$}(.4,1.2)
            edge[bend right]node[matrix on edge,pos=.65]{$A$}(.7,1.2)
            edge[bend right](-.7,-.5)
            edge[bend right](-.4,-.5)
            edge[bend left](.7,-.5);
        \draw(.05,1.2)node[toplabel]{$I_1^c\:\cdots\:I_{n-k}^c$};
        \draw(0,-.5)node[bottomlabel]{$J_1J_2\:\cdots\:J_k$};
    }
    =\sgn(J\:\overset\leftarrow{J^c})\frac{\sgn(I^c\:\overset\leftarrow{I})}{(n-k)!}
    \tikz[heighttwo,scale=1.2]{
        \node[ciliatednode=180](vb)at(0,1.6){}
            edge[bend left](-.7,2.2)
            edge[bend left](-.4,2.2)
            edge[bend right](.7,2.2);
        \node[ciliatednode=180](va)at(0,.3){}
            edge[out=160,in=200,looseness=2]node[matrix on edge]{$A$}(vb)
            edge[out=20,in=-20,looseness=2]node[matrix on edge]{$A$}(vb)
            edge[out=20,in=-20,looseness=1.2]node[matrix on edge]{$A$}(vb)
            edge[bend right](-.7,-.2)
            edge[bend right](-.4,-.2)
            edge[bend left](.7,-.2);
        \node[bottomlabel]at(0,-.2){$J_1J_2\cdots J_k$};
        \node[toplabel]at(0,2.2){$I_1I_2\cdots I_k$};
        \draw[dotdotdot](-.6,1)to node[toplabel,scale=.8]{$n\!-\!k$}(.4,1);
    }.
    \end{equation}
    Matching this up with the cofactor $C_{I,J}=(-1)^{I_1+\cdots+I_k+J_1+\cdots+J_k}[A^c_{I,J}]$ requires a little bit of work with the signs.
    \begin{lemma}
        Let $J=(J_1,\ldots,J_k)$ and $J^c=(J^c_1,\ldots,J^c_{n-k})$ be ordered increasing subsets of $N$ whose union is $N$. Then
            $$\mathrm{sgn}(J^c\:\overset\leftarrow{J}) = (-1)^{nk+J_1+J_2+\cdots+J_k}.$$
    \begin{proof}
        Move the $\{J_i\}$ one at a time to their ``proper'' positions among the $J^c$. The ordering implies
            $$(\ldots,J^c_{J_k-k+1},\ldots,J^c_{n-k},J_k,\ldots)=(\ldots,J_k+1,\ldots,n,J_k,\ldots),$$
        so $n-J_k$ transpositions are required to return $J_k$ to its proper place. Repeating this for each other $J_i$ gives the identity after a total of $nk-(J_1+\cdots+J_k)$ transpositions.
    \end{proof}
    \end{lemma}
    Thus $\mathrm{sgn}(J\:\overset\leftarrow{J^c})\mathrm{sgn}(I^c\:\overset\leftarrow{I})
    =\tdsign{n}(-1)^{I_1+\cdots+I_k+J_1+\cdots+J_k}$, verifying the diagram for the general cofactor is as stated.

    The adjugate diagram is the case $k=1$ with matrix orientations reversed to handle the transpose.
\end{proof}
\end{proposition}

\subsection{Decomposition of Trace Diagrams}

\begin{definition}\label{d:tdminor}
    Given a matrix $A$, a \emph{diagram $A$-minor} is an (unframed) diagram with a single $n$-vertex in which a subset of the edges may be labeled by $A$, in such a way that all matrix markings are compatibly oriented. In particular, the diagram may be written as $\pm1$ times a diagram of the form
    $$
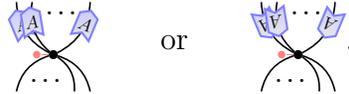
\tikz{
        \node[ciliatednode=180](va)at(0,.2){}
            edge[bend left]node[reverse matrix on edge,pos=.65]{$A$}(-.7,1.2)
            edge[bend left]node[reverse matrix on edge,pos=.65]{$A$}(-.4,1.2)
            edge[bend right]node[matrix on edge,pos=.65]{$A$}(.7,1.2)
            \foreach\xa/\xb in{right/-.7,left/.4,left/.7}{edge[bend \xa](\xb,-.5)};
        \draw[dotdotdot](-.4,1)to(.7,1);\draw[dotdotdot](-.7,-.3)to(.4,-.3);}
    \qquad\text{or}\qquad
    \tikz{
        \node[ciliatednode=180](va)at(0,.2){}
            edge[bend left]node[matrix on edge,pos=.65]{$A$}(-.7,1.2)
            edge[bend left]node[matrix on edge,pos=.65]{$A$}(-.4,1.2)
            edge[bend right]node[reverse matrix on edge,pos=.65]{$A$}(.7,1.2)
            \foreach\xa/\xb in{right/-.7,left/.4,left/.7}{edge[bend \xa](\xb,-.5)};
        \draw[dotdotdot](-.4,1)to(.7,1);\draw[dotdotdot](-.7,-.3)to(.4,-.3);}.
    $$
    (The sign comes from the possible need to switch the order of edges at the $n$-vertex so that all edges with matrices are adjacent.)
\end{definition}

Proposition \ref{p:minor-diagram} states that a diagram $A$-minor evaluates to a matrix minor $\pm[A_{I,J}]$ when the ends of the strands are labeled by $I$ and $J$. The following theorem states the conditions under which a trace diagram may be decomposed into diagram minors.

\begin{theorem}\label{thm:diagram-minors}
    Let $\mathcal{D}$ be a trace diagram in which every matrix marking is adjacent to an $n$-vertex. Then $\mathcal{D} = C \mathcal{D}'$ for some $\mathcal{D}'$ that may be decomposed into diagram minors, where $C$ is a constant that does not depend on any matrix entries.
\begin{proof}
    In this proof ``equivalence'' will mean equal up to a constant factor that does not depend on any matrix entries. The key step in the theorem is to use the cut-and-paste lemma to introduce additional $n$-vertices as necessary to separate matrices by node. For instance, the diagram
    $$
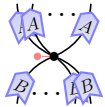
\tikz{
        \node[ciliatednode=180](va)at(0,.2){}
            edge[bend left]node[reverse matrix on edge,pos=.65]{$A$}(-.7,1.2)
            edge[bend left]node[reverse matrix on edge,pos=.65]{$A$}(-.4,1.2)
            edge[bend right]node[matrix on edge,pos=.65]{$A$}(.7,1.2)
            edge[bend right]node[matrix on edge,pos=.7]{$B$}(-.7,-.7)
            edge[bend left]node[reverse matrix on edge,pos=.7]{$B$}(.4,-.7)
            edge[bend left]node[reverse matrix on edge,pos=.7]{$B$}(.7,-.7);
        \draw[dotdotdot](-.4,1)to(.7,1);\draw[dotdotdot](-.7,-.5)to(.4,-.5);}
    $$
    cannot be decomposed into minors. However, using the cut-and-paste lemma and Proposition \ref{p:two-node-function}, is equivalent to
    $$
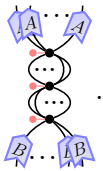
\tikz[heighttwo]{
        \node[ciliatednode=180](va)at(0,1.7){}
            edge[bend left]node[reverse matrix on edge,pos=.65]{$A$}(-.7,2.7)
            edge[bend left]node[reverse matrix on edge,pos=.65]{$A$}(-.4,2.7)
            edge[bend right]node[matrix on edge,pos=.65]{$A$}(.7,2.7);
        \node[ciliatednode=180](vc)at(0,1){};
        \node[ciliatednode=180](vb)at(0,.3){}
            edge[bend right]node[matrix on edge,pos=.7]{$B$}(-.7,-.7)
            edge[bend left]node[reverse matrix on edge,pos=.7]{$B$}(.4,-.7)
            edge[bend left]node[reverse matrix on edge,pos=.7]{$B$}(.7,-.7);
        \draw[dotdotdot](-.7,-.5)to(.4,-.5);\draw[dotdotdot](-.4,2.5)to(.7,2.5);
        \draw(va)to[out=200,in=160,looseness=2](vc);
        \draw(va)to[out=-20,in=20,looseness=1.2](vc);
        \draw(va)to[out=-20,in=20,looseness=2](vc);
        \draw(vc)to[out=200,in=160,looseness=2](vb);
        \draw(vc)to[out=200,in=160,looseness=1.22](vb);
        \draw(vc)to[out=-20,in=20,looseness=2](vb);
        \draw[dotdotdot](-.3,.65)to(.5,.65);\draw[dotdotdot](-.5,1.35)to(.3,1.35);
    }.
    $$
    Proceeding in this manner, since every matrix is adjacent to an $n$-vertex, one may introduce enough vertices in $\mathcal{D}$ to obtain an equivalent diagram $\mathcal{D}'$ such that every $n$-vertex in $\mathcal{D}$ is adjacent to a unique matrix with consistent orientation. One may then cut around each $n$-vertex in a diagram, including the adjacent matrices, to decompose the diagram into diagram minors.
\end{proof}
\end{theorem}

It follows immediately from this theorem that any such diagram may be expressed as a polynomial function of matrix minors. This in itself is not surprising, since the entries of a matrix are technically minors. The power of the result is that the structure of trace diagrams allows one to accomplish this decomposition ``efficiently'' by giving an upper bound for the number of minors in the decomposition.

For the purposes of the next theorem, we say that a collection of matrix markings form a \emph{compatible matrix collection} if (i) they have the same matrix label, (ii) they are adjacent to the same $n$-vertex, and (iii) they have the same orientation relative to the $n$-vertex. Given a trace diagram $\mathcal{D}$ in which every matrix is adjacent to an $n$-vertex, define the \emph{compatible partition number $N_{\mathcal{D}}$} of a trace diagram to be the minimum number of collections in a partition of all matrix markings in a diagram into compatible collections. For example,
    $$
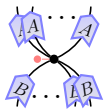
\tikz{
        \node[ciliatednode=180](va)at(0,.2){}
            edge[bend left]node[reverse matrix on edge,pos=.65]{$A$}(-.7,1.2)
            edge[bend left]node[reverse matrix on edge,pos=.65]{$A$}(-.4,1.2)
            edge[bend right]node[matrix on edge,pos=.65]{$A$}(.7,1.2)
            edge[bend right]node[matrix on edge,pos=.7]{$B$}(-.7,-.7)
            edge[bend left]node[reverse matrix on edge,pos=.7]{$B$}(.4,-.7)
            edge[bend left]node[reverse matrix on edge,pos=.7]{$B$}(.7,-.7);
        \draw[dotdotdot](-.4,1)to(.7,1);\draw[dotdotdot](-.7,-.5)to(.4,-.5);}
    $$
contains two compatible matrix collections, and the compatible partition number is 2.

\begin{theorem}\label{thm:diagram-minors-2}
    Let $\mathcal{D}$ be a trace diagram in which every matrix marking is adjacent to a vertex, and let $N_{\mathcal{D}}$ be the compatible partition number of $\mathcal{D}$. Then, the trace diagram function $f_{\mathcal{D}}$ may be expressed as a summation over a product of $N_{\mathcal{D}}$ matrix minors.
\begin{proof}
    In the proof of Theorem \ref{thm:diagram-minors}, one may ensure that every compatible matrix collection remains adjacent to the same vertex. Thus, one may write $\mathcal{D}=C\mathcal{D}'$, where $\mathcal{D}'$ decomposes into $N_{\mathcal{D}}$ diagram minors (and possibly some additional $n$-vertices without matrix markings). Given this decomposition, both $\mathcal{D}$ and $\mathcal{D}'$ may be expressed as summations over a product of $N_{\mathcal{D}}$ matrix minors.
\end{proof}
\end{theorem}

While $N_{\mathcal{D}}$ provides an upper bound for the minimum number of minors, it is not necessarily sharp. For example, the diagram
    $$
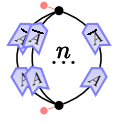
\tikz[heighttwo]{
        \node[ciliatednode=170](vb)at(0,2){};
        \node[ciliatednode=190](va)at(0,0){}
            edge[out=160,in=200,looseness=1.5]node[reverse matrix on edge,pos=.3,scale=.9]{$A$}node[matrix on edge,pos=.7,scale=.9]{$\overline{A}$}(vb)
            edge[out=150,in=210,looseness=1]node[reverse matrix on edge,pos=.3,scale=.9]{$A$}node[matrix on edge,pos=.7,scale=.9]{$\overline{A}$}(vb)
            edge[out=20,in=-20,looseness=1.5]node[matrix on edge,pos=.3,scale=.9]{$A$}node[reverse matrix on edge,pos=.7,scale=.9]{$\overline{A}$}(vb);
        \draw[dotdotdot](-.4,.9)to node[toplabel]{$n$}(.6,.9);
    }$$
has a compatible partition number of 2, but evaluates to $\tdsign{n} n!$.

\section{Three Short Determinant Proofs}\label{s:proofs}

There are several standard methods for computing the determinant. The \emph{Leibniz rule} is the common definition using permutations. \emph{Cofactor expansion} provides a recursive technique that lends itself well to by-hand calculations. \emph{Laplace expansion} is similar but uses generalized cofactors. A lesser known technique is \emph{Dodgson condensation} \cite{Dod1866}, which involves recursive computations using $2\times2$ determinants.

Diagrammatic techniques can unify these various approaches. Theorem \ref{thm:diagram-minors} leads to a straightforward diagrammatic approach to finding determinant identities: decompose the diagram for the determinant into pieces containing at most one node, and express the result as a summation over matrix minors. This approach gives the cofactor and Laplace formulae.

\subsection{Cofactor and Laplace Expansion}

\begin{proposition}[cofactor expansion]\label{p:cofactor-expansion}
    For an $n\times n$ matrix $A$ and $j\in\{1,2,\ldots,n\}$,
    \begin{equation}\label{eq:cofactor-expansion}
        \det(A) =\sum_{i=1}^n a_{ij} C_{ij} =\sum_{i=1}^n a_{ji} C_{ji}.
    \end{equation}
\begin{proof}
    Proposition \ref{p:determinant-diagram} states that
    $$\tikz[heightoneonehalf]{
            \node[ciliatednode=170](topnode)at(0,1.5){};
            \node[ciliatednode=-170](bottomnode)at(0,0){};
            \draw(bottomnode)to[bend left](-.7,.75)node[small matrix]{$A$}to[bend left](topnode);
            \draw(bottomnode)to[bend left](-.4,.75)node[small matrix]{$A$}to[bend left](topnode);
            \draw(bottomnode)to[bend right](.7,.75)node[small matrix]{$A$}to[bend right](topnode);
            \draw[dotdotdot](-.3,.75)to node[toplabel]{$n$}(.6,.75);
        }
    =(-1)^{\lfloor\frac{n}{2}\rfloor} n! \det(A).$$
    The diagram for the cofactor was found in Proposition \ref{p:cofactor-diagram}. The main idea in the proof is that it is possible to label one strand of the diagram arbitrarily, a consequence of two applications of the cut-and-paste lemma \eqref{eq:cut-paste}:
    $$\tikz[heightoneonehalf]{
            \node[ciliatednode=170](topnode)at(0,1.5){};
            \node[ciliatednode=-170](bottomnode)at(0,0){};
            \draw(bottomnode)to[bend left](-.7,.75)node[small matrix]{$A$}to[bend left](topnode);
            \draw(bottomnode)to[bend left](-.4,.75)node[small matrix]{$A$}to[bend left](topnode);
            \draw(bottomnode)to[bend right](.7,.75)node[small matrix]{$A$}to[bend right](topnode);
            \draw[dotdotdot](-.3,.75)to node[toplabel]{$n$}(.6,.75);
        }
    = n!\:\sgn(\beta)
    \tikz[scale=1.4,shift={(0,-.2)}]{
        \node[ciliatednode=180](va)at(0,0){}
            edge[bend left]node[reverse matrix on edge]{$A$}(-.7,1)
            edge[bend left]node[reverse matrix on edge]{$A$}(-.4,1)
            edge[bend right]node[matrix on edge]{$A$}(.7,1);
        \node[toplabel]at(0,1){$\beta_1\beta_2\:\:\cdots\:\beta_n$};
    } = \frac{n!}{(n-1)!} \sgn(\beta)^2
    \tikz[heighttwo,scale=1.2]{
        \node[ciliatednode=160](vb)at(0,1.8){}edge node[rightlabel,pos=1]{$\beta_n$}(0,2.2);
        \node[ciliatednode=200](va)at(0,.8){}
            edge[out=160,in=200,looseness=1.7]node[matrix on edge]{$A$}(vb)
            edge[out=20,in=-20,looseness=1.6]node[matrix on edge]{$A$}(vb)
            edge[out=20,in=-20,looseness=1.1]node[matrix on edge]{$A$}(vb)
            edge node[matrix on edge]{$A$}node[rightlabel,pos=1]{$\beta_n$}(0,-.2);
        \draw[dotdotdot](-.4,1.3)to(.2,1.3);
    }=n\: \tikz[heighttwo,scale=1.2]{
        \node[ciliatednode=160](vb)at(0,1.8){}edge node[rightlabel,pos=1]{$i$}(0,2.2);
        \node[ciliatednode=200](va)at(0,.8){}
            edge[out=160,in=200,looseness=1.7]node[matrix on edge]{$A$}(vb)
            edge[out=20,in=-20,looseness=1.6]node[matrix on edge]{$A$}(vb)
            edge[out=20,in=-20,looseness=1.1]node[matrix on edge]{$A$}(vb)
            edge node[matrix on edge]{$A$}node[rightlabel,pos=1]{$i$}(0,-.2);
        \draw[dotdotdot](-.4,1.3)to(.2,1.3);
    },$$    
    where $i=\beta_n$. This diagram may be evaluated by summing along an interior strand:
    $$n
    \tikz[heighttwo,scale=1.2]{
        \node[ciliatednode=160](vb)at(0,1.8){}edge node[rightlabel,pos=1]{$i$}(0,2.2);
        \node[ciliatednode=200](va)at(0,.8){}
            edge[out=160,in=200,looseness=1.7]node[matrix on edge]{$A$}(vb)
            edge[out=20,in=-20,looseness=1.6]node[matrix on edge]{$A$}(vb)
            edge[out=20,in=-20,looseness=1.1]node[matrix on edge]{$A$}(vb)
            edge node[matrix on edge]{$A$}node[rightlabel,pos=1]{$i$}(0,-.2);
        \draw[dotdotdot](-.4,1.3)to(.2,1.3);
    }
    =n\sum_{j=1}^n
    \tikz[heighttwo,scale=1.2]{
        \node[ciliatednode=160](vb)at(0,1.7){}edge node[rightlabel,pos=1]{$i$}(0,2.2);
        \node[ciliatednode=200](va)at(0,.5){}
            edge[out=160,in=200,looseness=1.7]node[matrix on edge]{$A$}(vb)
            edge[out=20,in=-20,looseness=1.6]node[matrix on edge]{$A$}(vb)
            edge[out=20,in=-20,looseness=1.1]node[matrix on edge]{$A$}(vb)
            edge node[rightlabel,pos=1]{$j$}(0,0);
        \draw[dotdotdot](-.4,1.1)to(.2,1.1);
    }
    \:\:
    \tikz{\draw(0,1.2)node[rightlabel]{$j$}to node[matrix on edge]{$A$}(0,-.2)node[rightlabel]{$i$};}
    = \tdsign{n} n(n-1)! \sum_{j=1}^n C_{ij}a_{ij}.
    $$
    Canceling the common $\tdsign{n}n!$ factor proves the first equality. The second equality follows by transposing the diagrams.
\end{proof}
\end{proposition}

This result is easily generalized by labeling several strands instead of just one (for a classical proof of this result, see Theorem 1 in Section 2.4 of \cite{LT85}).
\begin{proposition}[Laplace expansion]\label{p:laplace-expansion}
    $$\det(A)= \sum_{1\le J_1<\cdots<J_k\le n} C_{I,J} [A_{I,J}]= \sum_{1\le J_1<\cdots<J_k\le n} C_{J,I} [A_{J,I}].$$
\begin{proof}
    The proof is a variation of the above, this time cutting open the diagram along $k$ strands. First,
    \begin{equation*}
    \tikz[heightoneonehalf]{
            \node[ciliatednode=170](topnode)at(0,1.5){};
            \node[ciliatednode=-170](bottomnode)at(0,0){};
            \draw(bottomnode)to[bend left](-.7,.75)node[small matrix]{$A$}to[bend left](topnode);
            \draw(bottomnode)to[bend left](-.4,.75)node[small matrix]{$A$}to[bend left](topnode);
            \draw(bottomnode)to[bend right](.7,.75)node[small matrix]{$A$}to[bend right](topnode);
            \draw[dotdotdot](-.3,.75)to node[toplabel]{$n$}(.6,.75);
        }
    =\tfrac{n!}{(n-k)!}
    \tikz[heighttwo,scale=1.2]{
        \node[ciliatednode=180](vb)at(0,1.8){}
            edge[bend left](-.7,2.2)
            edge[bend left](-.4,2.2)
            edge[bend right](.7,2.2);
        \node[ciliatednode=180](va)at(0,.8){}
            edge[out=160,in=200,looseness=2]node[matrix on edge]{$A$}(vb)
            edge[out=20,in=-20,looseness=2]node[matrix on edge]{$A$}(vb)
            edge[out=20,in=-20,looseness=1.2]node[matrix on edge]{$A$}(vb)
            edge[bend right]node[reverse matrix on edge]{$A$}(-.7,-.2)
            edge[bend right]node[reverse matrix on edge]{$A$}(-.4,-.2)
            edge[bend left]node[matrix on edge]{$A$}(.7,-.2);
        \node[bottomlabel]at(0,-.2){$I_1I_2\cdots I_k$};
        \node[toplabel]at(0,2.2){$I_1I_2\cdots I_k$};
        \draw[dotdotdot](-.6,1.3)to(.4,1.3);
    }
    \end{equation*}
    We now use the cut-and-paste-lemma \eqref{l:cut-paste} to add an additional node at the bottom of the diagram, and then express the diagram as a summation over the interior labels to obtain:
    \begin{equation*}
    \tfrac{n!}{(n-k)!k!}
    \sgn(I^c\:\overset\leftarrow{I})
    \tikz[heightthree,scale=1.4]{
        \node[ciliatednode=180](vb)at(0,2.8){}
            edge[bend left](-.7,3.2)
            edge[bend left](-.4,3.2)
            edge[bend right](.7,3.2);
        \node[ciliatednode=180](va)at(0,1.5){}
            edge[out=160,in=200,looseness=2]node[matrix on edge]{$A$}(vb)
            edge[out=20,in=-20,looseness=2]node[matrix on edge]{$A$}(vb)
            edge[out=20,in=-20,looseness=1.2]node[matrix on edge]{$A$}(vb);
        \node[ciliatednode=180](vc)at(0,.2){}
            edge[out=160,in=200,looseness=2]node[reverse matrix on edge]{$A$}(va)
            edge[out=160,in=200,looseness=1.2]node[reverse matrix on edge]{$A$}(va)
            edge[out=20,in=-20,looseness=2]node[reverse matrix on edge]{$A$}(va)
            edge[bend right](-.7,-.2)
            edge[bend left](.4,-.2)
            edge[bend left](.7,-.2);
        \node[bottomlabel]at(0,-.2){$I_1^c\:\cdots\:I_k^c$};
        \node[toplabel]at(0,3.2){$I_1I_2\cdots I_k$};
        \draw[dotdotdot](-.6,2.1)to(.4,2.1);
        \draw[dotdotdot](-.4,.9)to(.6,.9);
    }
    =
    \tfrac{n!k!}{(n-k)!k!}
    \sgn(I^c\:\overset\leftarrow{I})
    \sum_{1\le J_1<\cdots<J_k\le n}
    \tikz[heighttwo,scale=1.2]{
        \node[ciliatednode=180](vb)at(0,1.6){}
            \foreach\xa/\xb in{left/-.7,right/.4,right/.7}{edge[bend \xa](\xb,2.2)};
        \node[ciliatednode=180](va)at(0,.3){}
            edge[out=160,in=200,looseness=2]node[matrix on edge]{$A$}(vb)
            edge[out=20,in=-20,looseness=2]node[matrix on edge]{$A$}(vb)
            edge[out=20,in=-20,looseness=1.2]node[matrix on edge]{$A$}(vb)
            edge[bend right](-.7,-.2)
            edge[bend right](-.4,-.2)
            edge[bend left](.7,-.2);
        \node[bottomlabel]at(0,-.2){$J_1J_2\cdots J_k$};
        \node[toplabel]at(0,2.2){$I_1I_2\cdots I_k$};
        \draw[dotdotdot](-.6,1)to node[toplabel,scale=.8]{$n\!-\!k$}(.4,1);
    }
    \tikz{
        \node[ciliatednode=180](va)at(0,.2){}
            edge[bend left]node[matrix on edge,pos=.65]{$A$}(-.7,1.2)
            edge[bend left]node[matrix on edge,pos=.65]{$A$}(-.4,1.2)
            edge[bend right]node[reverse matrix on edge,pos=.65]{$A$}(.7,1.2)
            edge[bend right](-.7,-.5)
            edge[bend left](.4,-.5)
            edge[bend left](.7,-.5);
        \draw(0,1.2)node[toplabel]{$J_1J_2\cdots J_k$};
        \draw(.2,-.5)node[bottomlabel]{$I_1^c\:\cdots\:I_{n-k}^c$};
    }.
    \end{equation*}
    By Propositions \ref{p:cofactor-diagram} and \ref{p:minor-diagram}, the first diagram here is $\tdsign{n}(n-k)! C_{I,J}$, and the second is $\sgn(I^c\:\overset\leftarrow{I}) [A_{I,J}]$. Matching up terms, we have now proven that  $\det(A)= \sum_{1\le J_1<\cdots<J_k\le n} C_{I,J} [A_{I,J}]$. The second statement is proven similarly.
\end{proof}
\end{proposition}

\subsection{A Determinant Theorem of Jacobi}

We now turn to the Jacobi determinant theorem, first stated in \cite{Jac1841}, which is used to derive Dodgson condensation \cite{RT07}. In contrast to the previous proofs, we state first the diagrammatic theorem, and show Jacobi's result as a corollary. This proof was first given in \cite{Mo08}.
\begin{proposition}\label{p:jacobi-diagrammatic}
    Let $A$ be an invertible $n\times n$ matrix, and let $I$ and $J$ be ordered subsets of $N$. Then
    \begin{equation}\label{dia:detadj}
    \begin{tikzpicture}[heighttwo,shift={(0,-2)}]
        \node[ciliatednode=-175](bottomnode)at(0,1){}
            \foreach\xa/\xb in{right/-.9,right/-.5,left/.9}{edge[bend \xa](\xb,0)};
        \node[ciliatednode=175](topnode)at(0,5){}
            \foreach\xa/\xb/\xc in{left/-.9,left/-.5,right/.9}{edge[bend \xa](\xb,6)};
        \foreach\xx/\xxb/\xxc in{-2.2/right/left,.8/left/right,2.6/left/right}{
            \node[ciliatednode=-150](tbottom)at(\xx,2.3){}edge[bend \xxb](bottomnode);
            \node[ciliatednode=150](ttop)at(\xx,3.7){}
                \foreach\xa/\xb in{left/1.5,left/.75,right/1.5}{%
                    edge[bend \xa=90,looseness=\xb]node[matrix on edge]{$A$}(tbottom)}
                edge[bend \xxc](topnode);
            \draw[dotdotdot,shift={(\xx,1.5)}](-.7,1.5)to node[basiclabel,above=2pt,scale=.8]{$n\!-\!1$}(.4,1.5);
        }
        \draw[dotdotdot](-1.7,4)to node[basiclabel,above,scale=.8]{$k$}(.3,4);
        \node[toplabel]at(0,6){$I_1^cI_2^c\cdots I_{n-k}^c$};
        \node[bottomlabel]at(0,0){$J_1^cJ_2^c\cdots J_{n-k}^c$};
    \end{tikzpicture}
    =c_1c_2 \det(A)^{k-1}
    \tikz[heighttwo,scale=1.2]{
        \node[ciliatednode=180](vb)at(0,1.6){}
            \foreach\xa/\xb in{left/-.7,left/-.4,right/.7}{edge[bend \xa](\xb,2.2)};
        \node[ciliatednode=180](va)at(0,.3){}
            \foreach\xa/\xb/\xc in{160/200/2,20/-20/1.2,20/-20/2}{
                edge[out=\xa,in=\xb,looseness=\xc]node[reverse matrix on edge]{$A$}(vb)}
            \foreach\xa/\xb in{right/-.7,right/-.4,left/.7}{edge[bend \xa](\xb,-.2)};
        \node[bottomlabel]at(0,-.2){$J_1J_2\cdots J_k$};
        \node[toplabel]at(0,2.2){$I_1I_2\cdots I_k$};
        \draw[dotdotdot](-.6,1)to node[toplabel,scale=.8]{$n\!-\!k$}(.4,1);
    },
    \end{equation}
    where $c_1c_2=\left(\tdsign{n} (n-1)!\right)^k
        \mathrm{sgn}(J^c\:\overset\leftarrow{J})\mathrm{sgn}(I^c\:\overset\leftarrow{I})\frac{k!}{(n-k)!}$.
\begin{proof}
    Use Proposition \ref{p:node-equivariance} to move each group of $n-1$ matrices in the lefthand diagram of \eqref{dia:detadj} onto a single edge labeled by $\overline{A}=A^{-1}$, then use Proposition \ref{p:two-node-function} with $k=1$ to eliminate the ``bubbles'' in the graph, as follows:
    $$
    \tikz[heighttwo,scale=1.2]{
        \node[ciliatednode=160](vb)at(0,1.6){}edge(0,2.2);
        \node[ciliatednode=200](va)at(0,.3){}
            edge[out=160,in=200,looseness=1.7]node[reverse matrix on edge]{$A$}(vb)
            edge[out=20,in=-20,looseness=1.6]node[reverse matrix on edge]{$A$}(vb)
            edge[out=20,in=-20,looseness=1.1]node[reverse matrix on edge]{$A$}(vb)
            edge(0,-.2);
        \draw[dotdotdot](-.4,.9)to node[toplabel,scale=.7]{$n\!-\!1$}(.2,.9);
    }
    =\det(A)
    \tikz[heighttwo,scale=1.2]{
        \node[ciliatednode=160](vb)at(0,1.2){}edge node[matrix on edge]{$\overline{A}$}(0,2.2);
        \node[ciliatednode=200](va)at(0,.2){}
            \foreach\xa/\xb/\xc in{160/200/1.7,20/-20/1.6,20/-20/1.1}{edge[out=\xa,in=\xb,looseness=\xc](vb)}
            edge(0,-.2);
        \draw[dotdotdot](-.4,.6)to node[toplabel,scale=.7]{$n\!-\!1$}(.2,.6);
    }
    =\det(A) \tdsign{n} (n-1)! \:
        \tikz{\draw(0,-.2)to node[matrix on edge]{$\overline{A}$}(0,1.2);}.
    $$
    This reduces the diagram to
    \begin{equation*}
    c_1 \det(A)^k
    \tikz[heighttwo,scale=1.2]{
        \node[ciliatednode=180](vb)at(0,1.6){}
            \foreach\xa/\xb in{left/-.7,right/.4,right/.7}{edge[bend \xa](\xb,2.2)};
        \node[ciliatednode=180](va)at(0,.3){}
            \foreach\xa/\xb/\xc in{160/200/2,160/200/1.2,20/-20/2}{
                edge[out=\xa,in=\xb,looseness=\xc]node[matrix on edge]{$\overline{A}$}(vb)}
            \foreach\xa/\xb in{right/-.7,left/.4,left/.7}{edge[bend \xa](\xb,-.2)};
        \node[bottomlabel]at(.1,-.2){$J_1^c\:\cdots\:J_{n-k}^c$};
        \node[toplabel]at(.1,2.2){$I_1^c\:\cdots\:I_{n-k}^c$};
        \draw[dotdotdot](-.4,1)to node[toplabel,scale=.9]{$k$}(.6,1);
    }
    = c_1 \det(A)^{k-1}
    \tikz[heighttwo,scale=1.2]{
        \node[ciliatednode=180](vb)at(0,1.3){}
                edge[bend left]node[matrix on edge,pos=.7]{$A$}(-.7,2.2)
                edge[bend right]node[reverse matrix on edge,pos=.7]{$A$}(.4,2.2)
                edge[bend right]node[reverse matrix on edge,pos=.7]{$A$}(.7,2.2);
        \node[ciliatednode=180](va)at(0,.2){}
            \foreach\xa/\xb/\xc in{160/200/2,160/200/1.2,20/-20/2}{edge[out=\xa,in=\xb,looseness=\xc](vb)}
            \foreach\xa/\xb in{right/-.7,left/.4,left/.7}{edge[bend \xa](\xb,-.2)};
        \node[bottomlabel]at(.1,-.2){$J_1^c\:\cdots\:J_{n-k}^c$};
        \node[toplabel]at(.1,2.2){$I_1^c\:\cdots\:I_{n-k}^c$};
        \draw[dotdotdot](-.4,.75)to node[toplabel,scale=.9]{$k$}(.6,.75);
    }
    = c_1 c_2 \det(A)^{k-1}
    \tikz[heighttwo,scale=1.2]{
        \node[ciliatednode=180](vb)at(0,1.6){}
            \foreach\xa/\xb in{left/-.7,left/-.4,right/.7}{edge[bend \xa](\xb,2.2)};
        \node[ciliatednode=180](va)at(0,.3){}
            \foreach\xa/\xb/\xc in{160/200/2,20/-20/1.2,20/-20/2}{
                edge[out=\xa,in=\xb,looseness=\xc]node[reverse matrix on edge]{$A$}(vb)}
            \foreach\xa/\xb in{right/-.7,right/-.4,left/.7}{edge[bend \xa](\xb,-.2)};
        \node[bottomlabel]at(0,-.2){$J_1J_2\cdots J_k$};
        \node[toplabel]at(0,2.2){$I_1I_2\cdots I_k$};
        \draw[dotdotdot](-.6,1)to node[toplabel,scale=.8]{$n\!-\!k$}(.4,1);
    }.
    \end{equation*}
    The second step is also a consequence of Proposition \ref{p:node-equivariance}. The third step uses the cut-and-paste lemma \eqref{eq:cut-paste} twice. The constants are $c_1=\left((-1)^{\lfloor\frac{n}{2}\rfloor}(n-1)!\right)^k$ and $c_2=\mathrm{sgn}(J^c\:\overset\leftarrow{J})\mathrm{sgn}(I^c\:\overset\leftarrow{I})\frac{k!}{(n-k)!}$.
\end{proof}
\end{proposition}

\begin{corollary}[Jacobi Determinant Theorem]\label{c:jacobi-theorem}
    Let $A$ be an $n\times n$ invertible matrix, and let $A_{I,J}$ be a $k\times k$ submatrix of $A$. Then
    \begin{align}\label{eq:jacobi}
        [\adj(A)_{I,J}]=C_{J,I}\det(A)^{k-1},
    \end{align}
    where $[\adj(A)_{I,J}]$ is the corresponding minor of the adjugate of $A$.
\begin{proof}
    Rewrite \eqref{dia:detadj} as $\mathcal{D}_1 = c_1 c_2 \det(A)^{k-1}\mathcal{D}_2$. By \eqref{eq:cofactor-diagram},
    \begin{equation}\label{eq:jacobi1}
        \mathcal{D}_2=\tdsign{n}(n-k)!C_{J,I} \equiv c_3 C_{J,I}.
    \end{equation}

    To see the meaning of $\mathcal{D}_1$, consider the following restatement of \eqref{eq:minor2}:
    $$[A_{I,J}]=\frac{\sgn(J^c\:\overset\leftarrow{J})\sgn(I\:\overset\leftarrow{I^c})}{k!}
    \tikz[heighttwo,scale=1.2]{
        \node[ciliatednode=180](vb)at(0,1.6){}
            \foreach\xa/\xb in{left/-.7,right/.4,right/.7}{edge[bend \xa](\xb,2.2)};
        \node[ciliatednode=180](va)at(0,.3){}
            \foreach\xa/\xb/\xc in{160/200/2,160/200/1.2,20/-20/2}{
                edge[out=\xa,in=\xb,looseness=\xc]node[matrix on edge]{$A$}(vb)}
            \foreach\xa/\xb in{right/-.7,left/.4,left/.7}{edge[bend \xa](\xb,-.2)};
        \node[bottomlabel]at(.1,-.2){$J_1^c\:\cdots\:J_{n-k}^c$};
        \node[toplabel]at(.1,2.2){$I_1^c\:\cdots\:I_{n-k}^c$};
        \draw[dotdotdot](-.4,1)to node[toplabel,scale=.9]{$k$}(.6,1);
    }.$$
    From this, one obtains a diagram for $[\mathsf{adj}(A)_{I,J}]$ by replacing each $A$ with the adjugate diagram \eqref{eq:cofactor-diagram}. The result is a multiple of $\mathcal{D}_1$:
    \begin{equation}\label{eq:jacobi2}
    [\mathsf{adj}(A)_{I,J}]
        = \frac{\sgn(J^c\:\overset\leftarrow{J})\sgn(I\:\overset\leftarrow{I^c})\left(\tdsign{n}\right)^k}
        {k!\left((n-1)!\right)^k} \mathcal{D} \equiv c_4 \mathcal{D}_1.
    \end{equation}
    Combining \eqref{dia:detadj}, \eqref{eq:jacobi1}, and \eqref{eq:jacobi2} gives
    $$[\mathsf{adj}(A)_{I,J}] = c_4 \mathcal{D}_1 = c_1 c_2 c_4 \det(A)^{k-1} \mathcal{D}_2 = c_1c_2c_3c_4\det(A)^{k-1} C_{J,I}.$$
    It is straightforward to check that $c_1c_2c_3c_4=1$.
\end{proof}
\end{corollary}

The first proofs of this theorem took several pages to complete, and required careful attention to indices and matrix elements. A modern proof is given in \cite{RT07} that also takes several pages, and relies on expressing the minor as the determinant of an $n\times n$ matrix derived from $A$. By contrast, the diagrammatic portion of the proof (Proposition \ref{p:jacobi-diagrammatic}) contains the essence of the result and was relatively easy. The more difficult part was showing that the diagrammatic relation corresponded to the correct algebraic statement.

Many identities in linear algebra are simply special cases of this theorem. For example, when $I=J=N$, then $[A^c_{I,J}]=1$ trivially and so
    $$\det(\adj(A))=\det(A)^{n-1}.$$

Charles Dodgson' \emph{condensation method} \cite{Dod1866} also depends on this result. The following example shows the condensation method at work on a $4\times4$ determinant.
    $$
	\begin{vmatrix}
	-2 & -1 & -1 & -4 \\
	-1 & -2 & -1 & -6 \\
	-1 & -1 &  2 &  4 \\
 	2 &  1 & -3 & -8 \\
	\end{vmatrix}
    \longrightarrow
    \begin{vmatrix}
 	3 & -1 &  2 \\
	-1 & -5 &  8 \\
 	1 &  1 & -4 \\
	\end{vmatrix}
    \longrightarrow
    \begin{vmatrix}8&-2\\-4&6\end{vmatrix}
    \longrightarrow
    -8,
    $$
where -8 is the determinant of the original matrix. Each step involves taking $2\times2$ determinants, making the process easy to do by hand. However, the technique fails for some matrices since it involves division.

The method relies on the particular case $I=J=\{1,n\}$. Then $[A^c_{I,J}]$ is the determinant of the interior entries, and $[\adj(A)_{I,J}]=C_{11}C_{nn}-C_{1n}C_{n1}$, where $C_{ij}$ is the cofactor, so \eqref{eq:jacobi} becomes
\begin{equation}\label{eq:condensationidentity}
    \det(A)=\frac{C_{11}C_{nn}-C_{1n}C_{n1}}{\det(\mathrm{int}(A))}.
\end{equation}
For $3\times 3$ matrices, this is precisely Dodgson's method. Larger determinants are computed using several iterations of this formula. 


\section{Generalizations using Trace Diagrams}\label{s:generalization}

One of the advantages of using trace diagrams is the ease with which certain proofs are generalized. This is because, in contrast to traditional proofs, patterns in trace diagram proofs are more easily recognized. For example, the proof of Proposition \ref{p:jacobi-diagrammatic} is readily generalized when $ik\le n$ to the following:
\begin{proposition}\label{p:jacobi-diagrammatic-2}
    Let $A$ be an invertible $n\times n$ matrix, and let $I$ and $J$ be ordered subsets of $N$. Then
    \begin{equation}\label{dia:detadj-gen}
    \begin{tikzpicture}[heighttwo,shift={(0,-1)}]
        \node[ciliatednode=-175](bottomnode)at(0,1){}
            \foreach\xa/\xb in{right/-.9,right/-.5,left/.9}{edge[bend \xa](\xb,0)};
        \foreach\xx/\xxb/\xxc in{-2.2/right/left,.8/left/right,2.6/left/right}{
            \node[ciliatednode=-165](tbottom)at(\xx,2.3){}
                edge[bend \xxb,looseness=.25](bottomnode)
                edge[dotdotdot,bend \xxb,looseness=.75](bottomnode)
                edge[bend \xxb,looseness=1.5](bottomnode);
            \node[ciliatednode=170](ttop)at(\xx,3.7){}
                \foreach\xa/\xb in{left/1.5,left/.75,right/1.5}{%
                    edge[bend \xa=90,looseness=\xb]node[matrix on edge]{$A$}(tbottom)};
            \draw[dotdotdot,shift={(\xx,1.5)}](-.7,1.5)to node[basiclabel,above=2pt,scale=.8]{$n\!-\!i$}(.4,1.5);
        }
        \foreach\xa in{-2.2,.8,2.6}{
            \draw[shift={(\xa,0)}]
                (0,3.7)to[bend left](-.5,4.5)
                (0,3.7)to[bend left](-.2,4.5)
                (0,3.7)to[bend right](.5,4.5);
        }
        \draw[dotdotdot](-1.7,3.8)to node[basiclabel,above,scale=.8]{$k$}(.3,3.8);
        \node[toplabel]at(.3,4.5){$I_1\cdots I_i\quad\cdots\quad\cdots\qquad\cdots I_{ik}$};
        \node[bottomlabel]at(.2,0){$J_1^cJ_2^c\cdots J_{n-ik}^c$};
    \end{tikzpicture}
    = c_1c_2 \det(A)^{k-1}
    \tikz[heighttwo,scale=1.2]{
        \node[ciliatednode=180](vb)at(0,1.6){}
            \foreach\xa/\xb in{left/-.7,left/-.4,right/.7}{edge[bend \xa](\xb,2.2)};
        \node[ciliatednode=180](va)at(0,.3){}
            \foreach\xa/\xb/\xc in{160/200/2,20/-20/1.2,20/-20/2}{
                edge[out=\xa,in=\xb,looseness=\xc]node[reverse matrix on edge]{$A$}(vb)}
            \foreach\xa/\xb in{right/-.7,right/-.4,left/.7}{edge[bend \xa](\xb,-.2)};
        \node[bottomlabel]at(0,-.2){$J_1J_2\cdots J_{ik}$};
        \node[toplabel]at(0,2.2){$I_1I_2\cdots I_{ik}$};
        \draw[dotdotdot](-.6,1)to node[toplabel,scale=.7]{$n\!-\!ik$}(.4,1);
    },
    \end{equation}
    where $c_1c_2=\left(\tdsign{n} (n-i)!\right)^k\frac{\mathrm{sgn}(J\:\overleftarrow{J^c})}{(n-ik)!}$.
\begin{proof}
    The proof is similar to that of Proposition \ref{p:jacobi-diagrammatic}. Begin by reducing the diagram at left by applying the following steps at each small collection of $n-i$ matrices in the diagram:
    $$
    \tikz[heighttwo,scale=1.2]{
        \node[ciliatednode=180](vb)at(0,1.6){}
            \foreach\xa/\xb in{left/-.7,left/-.4,right/.7}{edge[bend \xa](\xb,2.2)};
        \node[ciliatednode=180](va)at(0,.3){}
            \foreach\xa/\xb/\xc in{160/200/2,20/-20/1.2,20/-20/2}{
                edge[out=\xa,in=\xb,looseness=\xc]node[reverse matrix on edge]{$A$}(vb)}
            \foreach\xa/\xb in{right/-.7,right/-.4,left/.7}{edge[bend \xa](\xb,-.2)};
        \draw[dotdotdot](-.6,1)to node[toplabel,scale=.7]{$n\!-\!i$}(.4,1);
        \draw[dotdotdot](-.4,2.2)to node[toplabel]{$i$}(.6,2.2);
        \draw[dotdotdot](-.4,-.2)to node[bottomlabel]{$i$}(.6,-.2);
    }
    = \det(A)
    \tikz[heighttwo,scale=1.2]{
        \node[ciliatednode=180](vb)at(0,1.2){}
            edge[bend left]node[reverse matrix on edge]{$\overline{A}$}(-.7,2.2)
            edge[bend left]node[reverse matrix on edge]{$\overline{A}$}(-.4,2.2)
            edge[bend right]node[matrix on edge]{$\overline{A}$}(.7,2.2);
        \node[ciliatednode=180](va)at(0,.2){}
            \foreach\xa/\xb/\xc in{160/200/2,20/-20/1.2,20/-20/2}{edge[out=\xa,in=\xb,looseness=\xc](vb)}
            \foreach\xa/\xb in{right/-.7,right/-.4,left/.7}{edge[bend \xa](\xb,-.2)};
        \draw[dotdotdot](-.6,.7)to node[toplabel,scale=.7]{$n\!-\!i$}(.4,.7);
        \draw[dotdotdot](-.4,2.2)to node[toplabel]{$i$}(.6,2.2);
        \draw[dotdotdot](-.4,-.2)to node[bottomlabel]{$i$}(.6,-.2);
    }
    \longrightarrow \det(A) \tdsign{n} (n-i)!\:
    \tikz[heighttwo]{
        \foreach\xa in{-.6,-.3, .6}{\draw(\xa,.2)to node[matrix on edge]{$\overline{A}$}(\xa,1.8);}
        \draw[dotdotdot](-.3,1)to node[toplabel]{$i$}(.6,1);}
    .$$
    Note that the last step is only true in the context of the larger diagram, in which case it follows by two applications of the cut-and-paste-lemma \eqref{eq:cut-paste}. After this step, the diagram reduces to
    \begin{equation*}
    c_1 \det(A)^k
    \tikz{
        \node[ciliatednode=180](va)at(0,.2){}
            edge[bend left]node[reverse matrix on edge,pos=.65]{$\overline{A}$}(-.7,1.2)
            edge[bend left]node[reverse matrix on edge,pos=.65]{$\overline{A}$}(-.4,1.2)
            edge[bend right]node[matrix on edge,pos=.65]{$\overline{A}$}(.7,1.2)
            \foreach\xa/\xb in{right/-.7,left/.4,left/.7}{edge[bend \xa](\xb,-.5)};
        \draw(0,1.2)node[toplabel]{$I_1I_2\cdots I_{ik}$};
        \draw(.2,-.5)node[bottomlabel]{$J_1^c\:\cdots\:J_{n-ik}^c$};
    }
    = c_1 c_2 \det(A)^{k-1}
     \tikz[heighttwo,scale=1.2]{
        \node[ciliatednode=180](vb)at(0,1.6){}
            \foreach\xa/\xb in{left/-.7,left/-.4,right/.7}{edge[bend \xa](\xb,2.2)};
        \node[ciliatednode=180](va)at(0,.3){}
            \foreach\xa/\xb/\xc in{160/200/2,20/-20/1.2,20/-20/2}{
                edge[out=\xa,in=\xb,looseness=\xc]node[reverse matrix on edge]{$A$}(vb)}
            \foreach\xa/\xb in{right/-.7,right/-.4,left/.7}{edge[bend \xa](\xb,-.2)};
        \node[bottomlabel]at(0,-.2){$J_1J_2\cdots J_{ik}$};
        \node[toplabel]at(0,2.2){$I_1I_2\cdots I_{ik}$};
        \draw[dotdotdot](-.6,1)to node[toplabel,scale=.7]{$n\!-\!ik$}(.4,1);
    },
    \end{equation*}
    where $c_1=\left(\tdsign{n}(n-i)!\right)^k$ and
    $c_2=\frac{\mathrm{sgn}(J\:\overleftarrow{J^c})}{(n-ik)!}$. The details here are identical to those in the proof of Proposition \ref{p:jacobi-diagrammatic}.
\end{proof}
\end{proposition}

We will use this result to prove a generalization of the Jacobi determinant theorem, which concerns a more general notion of a matrix minor. We must first introduce some new concepts. Let $V$ be an $n$-dimensional vector space. Given a multilinear transformation $\mathbf{A}:V^{\otimes i}\to V^{\otimes i}$, one can represent the value of the transformation by the coefficients
    $$(\mathbf{A})_{\alpha,\beta} \equiv \langle\bse_\alpha,\mathbf{A}\bse_\beta\rangle,$$
where $\alpha,\beta\in N^i$. Diagrammatically, $\mathbf{A}$ is represented by an oriented node with $i$ inputs and $i$ outputs:
    \tikz{\foreach\xa in{-.4,-.2,.4}{\draw(\xa,-.2)--(0,.5)--(\xa,1.2);}
        \draw[dotdotdot](-.1,0)to node[bottomlabel]{$i$}(.3,0);
        \node[small matrix,xscale=1.5]at(0,.5){$\mathbf{A}$};}.

The \emph{$i$-adjugate} of a matrix $A$ ($0\le i\le n$) is the multilinear transformation  $\adj_i(A):V^{\otimes i}\to V^{\otimes i}$
whose coefficients are general cofactors:
    $$\left(\adj_i(A)\right)_{I,J} = C_{J,I},$$
where $I$ and $J$ are ordered subsets of $N$ with $i$ elements. It follows from Proposition \ref{p:cofactor-diagram} that
    \begin{equation}\label{dia:general-adjugate}
    \adj_i(A)=
    \frac{\tdsign{n}}{(n-i)!}
    \tikz[heighttwo,scale=1.2]{
        \node[ciliatednode=180](vb)at(0,1.6){}
            \foreach\xa/\xb in{left/-.7,left/-.4,right/.7}{edge[bend \xa](\xb,2.2)};
        \node[ciliatednode=180](va)at(0,.3){}
            \foreach\xa/\xb/\xc in{160/200/2,20/-20/1.2,20/-20/2}{
                edge[out=\xa,in=\xb,looseness=\xc]node[reverse matrix on edge]{$A$}(vb)}
            \foreach\xa/\xb in{right/-.7,right/-.4,left/.7}{edge[bend \xa](\xb,-.2)};
        \draw[dotdotdot](-.6,1)to node[toplabel,scale=.7]{$n\!-\!i$}(.4,1);
        \draw[dotdotdot](-.4,2.2)to node[toplabel]{$i$}(.6,2.2);
        \draw[dotdotdot](-.4,-.2)to node[bottomlabel]{$i$}(.6,-.2);
    }.\end{equation}

We also need to generalize the idea of a matrix minor. Let $\mathbf{A}$ be a multilinear transformation, as defined above. Let a positive integer $k$ be chosen for which $0\le ik\le n$. Let $\mathbf{I}=(\mathbf{I}_1,\ldots,\mathbf{I}_k)$ consist of $k$ $i$-tuples with $\mathbf{I}_j\equiv(I_{j,1},\ldots,I_{j,i})$ and all elements of $\mathbf{I}$ distinct. Let the order of indices be chosen so that
    $$1\le I_{1,1}\le\cdots\le I_{1,i}\le\:\cdots\:\le I_{k,1}\le\cdots\le I_{k,i}.$$
Let $\mathbf{J}$ be similarly chosen. The \emph{$\mathbf{I},\mathbf{J}$-minor} of $\mathbf{A}$ is defined to be
    $$
    [\mathbf{A}_{\mathbf{I},\mathbf{J}}]
    =\sum_{\sigma\in S_{ik}} \sgn(\sigma)
        (\mathbf{A})_{\mathbf{I}_1,\sigma(\mathbf{J}_1)} (\mathbf{A})_{\mathbf{I}_2,\sigma(\mathbf{J}_2)}
        \cdots
        (\mathbf{A})_{\mathbf{I}_k,\sigma(\mathbf{J}_k)}.$$
Generalizing Proposition \ref{p:minor-diagram} gives
    \begin{equation}\label{dia:general-minor}
    [\mathbf{A}_{\mathbf{I},\mathbf{J}}]
        =\sgn(\mathbf{J}^c\:\overset\leftarrow{\mathbf{J}})
    \tikz{
        \node[ciliatednode=190](va)at(0,.2){}
            \foreach\xa/\xb in{right/-.7,right/-.4,left/.7}{edge[bend \xa](\xb,-.5)}
            \foreach\xa/\xb in{left/-1.8,left/-1.6,left/-1,left/-.6,left/-.4,right/.2,right/1.2,right/1.4,right/2}{
                edge[bend \xa](\xb,1)};
        \foreach\xsh in{-1.4,-.2,1.6}{
            \foreach\xa in{-.4,-.2,.4}{\draw[shift={(\xsh,1.2)}](\xa,-.2)--(0,.5)--(\xa,1.2);}
            \draw[dotdotdot,shift={(\xsh,1.2)}](-.1,0)to node[bottomlabel]{$i$}(.3,0);
            \node[small matrix,xscale=1.5]at(\xsh,1.7){$A$};
            }
        \draw(0,2.4)node[toplabel]{$I_1\cdots I_i\quad\cdots\quad\:\:\cdots I_{ik}$};
        \draw(.2,-.5)node[bottomlabel]{$J_1^c\:\cdots\:J_{n-ik}^c$};
    }.
    \end{equation}

We can now use the diagrammatic result \eqref{dia:detadj-gen} to generalize the Jacobi determinant theorem.
\begin{theorem}\label{thm:generalized-jacobi}
    Let $A$ be an $n\times n$ invertible matrix, and let $A_{\mathbf{I},\mathbf{J}}$ be an $ik\times ik$ submatrix of $A$. Then
    \begin{equation}\label{jacobi-gen}
        [\adj_i(A)_{\mathbf{I},\mathbf{J}}]=C_{\mathbf{J},\mathbf{I}}\det(A)^{k-1}.
    \end{equation}
\begin{proof}
    Rewrite \eqref{dia:detadj-gen} as $\mathcal{D}_1 = c_1 c_2 \det(A)^{k-1}\mathcal{D}_2$.
    As in the proof of the Jacobi determinant theorem (Corollary \ref{c:jacobi-theorem}), $\mathcal{D}_2=\tdsign{n}(n-ik)!C_{\mathbf{J},\mathbf{I}}\equiv c_3 C_{\mathbf{J},\mathbf{I}}$. The diagram $\mathcal{D}_1$ is obtained by inserting $k$ copies of the $i$-adjugate diagram \eqref{dia:general-adjugate} into the generalized minor diagram \eqref{dia:general-minor}, and so
        $$[\adj(A)_{\mathbf{I},\mathbf{J}}]
        = \left(\frac{\tdsign{n}}{(n-i)!}\right)^k\sgn(\mathbf{J}^c\overleftarrow{\mathbf{J}}) \mathcal{D}_1 \equiv c_4 \mathcal{D}_1.$$
    Combining these results, one has $[\adj(A)_{\mathbf{I},\mathbf{J}}]\equiv c_1c_2c_3c_4\det(A)^{k-1}C_{\mathbf{J},\mathbf{I}}$, and it is straightforward to verify that $c_1c_2c_3c_4=1$.
\end{proof}
\end{theorem}

\section{Final Remarks}\label{s:conclusion}

The main purpose of this paper has been to introduce the ideas of signed graph colorings and trace diagrams. A secondary purpose has been to provide a lexicon for their translation into linear algebra. The advantage in this approach to linear algebra lies in the ability to \emph{generalize} results, as was done in Section \ref{s:generalization}.

There is much more to be said about trace diagrams. The case $n=2$ was the starting point of the theory \cite{Lev56} and has been studied extensively, most notably providing the basis for spin networks \cite{CFS95,Kau91} and the Kauffman bracket skein module \cite{BFK96}. In the general case, the coefficients of the characteristic equation of a matrix can be understood as the $n+1$ ``simplest'' closed trace diagrams \cite{Pet08}.

The diagrammatic language also proves to be extremely useful in invariant theory. It allows for easy expression of the ``linearization'' of the characteristic equation \cite{Pet08}, from which several classical results of invariant theory are derived \cite{Dre07}. Diagrams have already given new insights in the theory of character varieties and invariant theory \cite{Bu97,LP09,Sik01}, and it is likely that more will follow.

\section*{Acknowledgments}

The authors would like to thank the referee for many valuable comments and suggestions for improvements. We also thank Paul Falcone, Bill Goldman, Sean Lawton, and James Lee for valuable discussions and Amanda Beecher, Jennifer Peterson, and Brian Winkel for their comments on early drafts of this paper.

%


\bibliographystyle{plain}
\bibliography{determinant}

\end{document}